\documentclass[12pt]{amsart}
\usepackage{amsmath, amssymb, latexsym, amsthm}
\usepackage{mathrsfs}
\usepackage[all]{xy}

      \newenvironment{changemargin}[2]{\begin{list}{}{
         \setlength{\topsep}{0pt}\setlength{\leftmargin}{0pt}
         \setlength{\rightmargin}{0pt}
         \setlength{\listparindent}{\parindent}
         \setlength{\itemindent}{\parindent}
         \setlength{\parsep}{0pt plus 1pt}
         \addtolength{\leftmargin}{#1}\addtolength{\rightmargin}{#2}
         }\item }{\end{list}}

\newcommand{\alephes}{{\aleph_0}}

\newcommand{\cl}[1]{\overline{#1}}
\newcommand{\Next}[1]{{#1^+}}
\newcommand{\intvl}[2]{{[#1(#2),#1(#2\!+\!1))}}

\newcommand{\PN}{{P(\N)}}
\newcommand{\Cantor}{{\{0,1\}^\N}}
\newcommand{\setseq}[1]{\{#1 : n\in\N\}}
\newcommand{\CH}{the Continuum Hypothesis}

\newcommand{\sr}[2]{{\txt{{${#1}$}\\{\tiny ${#2}$}}}}
\newcommand{\Dfin}{\mathfrak{D}_{\textrm{\rm fin}}}
\newcommand{\fD}{\mathfrak{D}}
\newcommand{\scrA}{\mathscr{A}}
\newcommand{\scrB}{\mathscr{B}}
\newcommand{\scrC}{\mathscr{C}}
\newcommand{\seq}[1]{\{#1\}_{n\in\N}}

\newcommand{\op}{\operatorname}
\newcommand{\maxfin}{\op{maxfin}}

\newcommand{\BG}{\B_\Gamma}
\newcommand{\BT}{\B_\mathrm{T}}
\newcommand{\BO}{\B_\Omega}
\newcommand{\BL}{\B_\Lambda}

\newcommand{\BOfat}{\B_{\Omega}^{\mathrm{fat}}}
\newcommand{\cD}{\mathcal{D}}
\newcommand{\cI}{\mathcal{I}}
\newcommand{\cW}{\mathcal{W}}

\newcommand{\cB}{\mathcal{B}}
\newcommand{\cU}{\mathcal{U}}
\newcommand{\cV}{\mathcal{V}}

\newcommand{\Tau}{\mathrm{T}}

\newcommand{\cF}{\mathcal{F}}
\newcommand{\cG}{\mathcal{G}}
\newcommand{\cM}{\mathcal{M}}
\newcommand{\cN}{\mathcal{N}}
\newcommand{\N}{\mathbb{N}}
\newcommand{\bbR}{\mathbb{R}}
\newcommand{\bbQ}{\mathbb{Q}}
\newcommand{\Null}{\mathcal{N}}
\newcommand{\NN}{{\N^\N}}
\newcommand{\NZ}{{\Z^\N}}
\newcommand{\NNup}{{\N^{\uparrow\N}}}
\newcommand{\roth}{{[\N]^{\aleph_0}}}
\newcommand{\cO}{\mathcal{O}}
\newcommand{\B}{\mathcal{B}}
\newcommand{\Q}{\rationals}
\newcommand{\R}{\reals}

\newcommand{\fu}{\mathfrak{u}}
\newcommand{\fr}{\mathfrak{r}}
\newcommand{\fod}{\mathfrak{o\!d}}
\newcommand{\Union}{\bigcup}
\newcommand{\Z}{\mathbb{Z}}
\newcommand{\Impl}{\Rightarrow}

\long\def\forget#1\forgotten{}

\newcommand{\fb}{\mathfrak{b}}
\newcommand{\fc}{\mathfrak{c}}
\newcommand{\fd}{\mathfrak{d}}

\newcommand{\fg}{\mathfrak{g}}
\newcommand{\itm}{\item}
\newcommand{\oo}{\infty}
\newcommand{\fp}{\mathfrak{p}}
\newcommand{\fs}{\mathfrak{s}}
\newcommand{\w}{\omega}
\newcommand{\x}{\times}

\newcommand{\Iff}{\Leftrightarrow}
\newcommand\comp{^{\text{\tt c}}}
\newcommand{\nin}{\notin}

\newcommand{\sbst}{\subseteq}
\newcommand{\spst}{\supseteq}
\newcommand{\sm}{\setminus}

\newcommand{\as}{\subseteq^*}

\newcommand{\E}{\exists}
\newcommand{\I}{\mathcal{I}}
\newcommand{\cJ}{\mathcal{J}}

\newcommand{\cov}{{\sf cov}}
\newcommand{\add}{{\sf add}}

\newcommand{\cf}{{\sf cf}}
\newcommand{\non}{{\sf non}}

\newcommand{\ft}{\mathfrak{t}}
\newcommand{\fh}{\mathfrak{h}}

\newtheorem{thm}{Theorem}[section]
\newtheorem{prop}[thm]{Proposition}
\newtheorem{prob}[thm]{Problem}
\newtheorem{lem}[thm]{Lemma}
\newtheorem{cor}[thm]{Corollary}

\theoremstyle{definition}
\newtheorem{defn}[thm]{Definition}

\theoremstyle{remark}
\newtheorem{rem}[thm]{Remark}

\newcommand{\be}{\begin{enumerate}}
\newcommand{\ee}{\end{enumerate}}
\newcommand{\bi}{\begin{itemize}}
\newcommand{\ei}{\end{itemize}}

\newcommand{\sone}{{\sf S}_1}    \newcommand{\sfin}{{\sf S}_{fin}}
\newcommand{\ufin}{{\sf U}_{fin}}
    
\newcommand{\Split}{\mathsf{Split}}

\newcommand{\reals}{{\mathbb R}}
\newcommand{\rationals}{{\mathbb Q}}

\author{Boaz Tsaban}
\thanks{Supported by the Koshland Center for Basic Research.}
\address{Department of Applied Mathematics and Computer Science,
Weizmann Institute of Science, Rehovot 76100, Israel}
\email{boaz.tsaban@weizmann.ac.il}
\urladdr{http://www.cs.biu.ac.il/\~{}tsaban}

\title[Additivity of covering properties]{Additivity numbers of covering properties}

\subjclass{%
Primary: 37F20; 
Secondary 26A03, 
03E75 
}

\keywords{%
Menger property, Hurewicz property, selection principles, additivity numbers,
Rudin-Keisler ordering, near coherence of filters%
}

\begin{document}
\begin{abstract}
The \emph{additivity number} of a topological property (relative to a given space)
is the minimal number of subspaces with this property whose union does not have the property.
The most well-known case is where this number is greater than $\aleph_0$,
i.e.\ the property is $\sigma$-additive.
We give a rather complete survey of the known results about the additivity
numbers of a variety of topological covering properties, including
those appearing in the Scheepers diagram (which contains, among others,
the classical properties of Menger, Hurewicz, Rothberger, and Gerlits-Nagy).
Some of the results proved here were not published beforehand, and many
open problems are posed.
\end{abstract}

\maketitle

\section{Introduction}

Assume that $\cI$ is a topological property.
For a topological space $X$, let $\cI(X)$
denote the subspaces of $X$ which possess the property $\cI$,
and assume that $\cup\cI(X)\nin\cI(X)$.
Define the \emph{additivity number of $\cI$} (relative to $X$) as
$$\add_X(\cI)=\min\{|\cF| : \cF\sbst\cI(X)\mbox{ and }\cup\cF\nin\cI(X)\}.$$
$\cI(X)$ is \emph{additive} when $\add_X(\cI)\ge\aleph_0$ and
\emph{$\sigma$-additive} when $\add_X(\cI)>\aleph_0$.
Sometimes it is useful to have more precise estimations of the additivity number
of a property, or even better, determine it exactly in terms of well-studied cardinals.
This is the purpose of this paper. We do that for a variety of topological covering
properties,
but some restriction is necessary.
We concentrate on the case that $X$ is separable, metrizable, and zero-dimensional.
This restriction allows for a convenient application of the combinatorial method.
Having established the results for this case, one can seek for generalizations
(which are sometimes straightforward).
Each topological space as above is homeomorphic to a set of irrational
numbers. Thus, it suffices to study $\add_{\bbR\sm\bbQ}(\cI)$, and we can therefore
omit the subscript.

\subsection{Covering properties}
Fix a space $X$.
An open cover $\cU$ of $X$ is
\emph{large} if each member of $X$ is contained in infinitely
many members of $\cU$.
$\cU$ is an \emph{$\omega$-cover} if $X\nin\cU$ and for
each finite $F\sbst X$, there is $U\in\cU$ such that $F\subseteq U$.
$\cU$ is a \emph{$\gamma$-cover} of $X$ if it is infinite and for each $x\in X$,
$x$ is a member of all but finitely many members of $\cU$.

Let $\cO$, $\Lambda$, $\Omega$, and $\Gamma$ denote the collections of all countable open
covers, large covers, $\omega$-covers, and $\gamma$-covers of $X$, respectively.
Similarly, let $\cB$, $\BL$, $\BO$, and $\BG$ denote the corresponding countable \emph{Borel}
covers of $X$.\footnote{By \emph{open cover} (respectively, \emph{Borel cover})
we mean a cover whose elements are open (respectively, Borel).}
Let
$\scrA$ and $\scrB$ be any of these classes. We consider the following three
properties which $X$ may or may not have.
\begin{itemize}
\item[$\sone(\scrA,\scrB)$:]
For each sequence $\seq{\cU_n}$ of members of $\scrA$,
there exist members $U_n\in\cU_n$, $n\in\N$, such that $\setseq{U_n}\in\scrB$.
\item[$\sfin(\scrA,\scrB)$:]
For each sequence $\seq{\cU_n}$
of members of $\scrA$, there exist finite
subsets $\cF_n\sbst\cU_n$, $n\in\N$, such that $\Union_{n\in\N}\cF_n\in\scrB$.
\item[$\ufin(\scrA,\scrB)$:]
For each sequence $\seq{\cU_n}$ of members of $\scrA$
which do not contain a finite subcover,
there exist finite subsets $\cF_n\sbst\cU_n$, $n\in\N$,
such that $\setseq{\cup\cF_n}\in\scrB$.
\end{itemize}

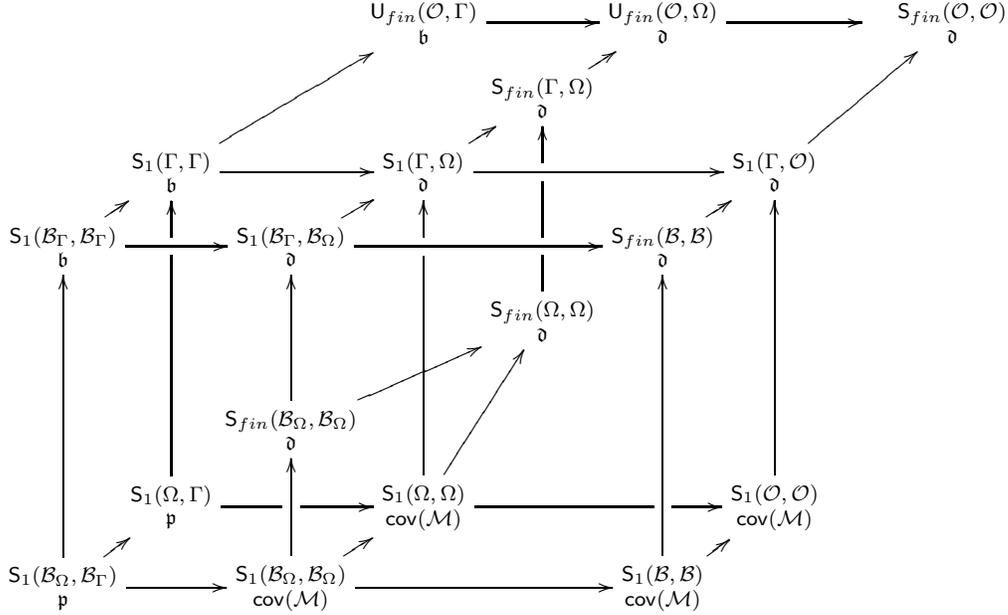
\begin{figure}[!htp]
\begin{changemargin}{-4cm}{-3cm}
\begin{center}
{\tiny
$\xymatrix@C=-2pt@R=6pt{
&
&
& \sr{\ufin(\cO,\Gamma)}{\fb}\ar[rr]
&
& \sr{\ufin(\cO,\Omega)}{\fd}\ar[rrrrr]
&
&
&
&
&
&
& \sr{\sfin(\cO,\cO)}{\fd}
\\
&
&
&
& \sr{\sfin(\Gamma,\Omega)}{\fd}\ar[ur]
\\
& \sr{\sone(\Gamma,\Gamma)}{\fb}\ar[rr]\ar[uurr]
&
& \sr{\sone(\Gamma,\Omega)}{\fd}\ar[rrr]\ar[ur]
&
&
& \sr{\sone(\Gamma,\cO)}{\fd}\ar[uurrrrrr]
\\
  \sr{\sone(\BG,\BG)}{\fb}\ar[ur]\ar[rr]
&
& \sr{\sone(\BG,\BO)}{\fd}\ar[ur]\ar[rrr]
&
&
& \sr{\sfin(\B,\B)}{\fd}\ar[ur]
\\
&
&
&
& \sr{\sfin(\Omega,\Omega)}{\fd}\ar'[u]'[uu][uuu]
\\
\\
&
& \sr{\sfin(\BO,\BO)}{\fd}\ar[uuu]\ar[uurr]
\\
& \sr{\sone(\Omega,\Gamma)}{\fp}\ar'[r][rr]\ar'[uuuu][uuuuu]
&
& \sr{\sone(\Omega,\Omega)}{\cov(\cM)}\ar'[uuuu][uuuuu]\ar'[rr][rrr]\ar[uuur]
&
&
& \sr{\sone(\cO,\cO)}{\cov(\cM)}\ar[uuuuu]
\\
  \sr{\sone(\BO,\BG)}{\fp}\ar[uuuuu]\ar[rr]\ar[ur]
&
& \sr{\sone(\BO,\BO)}{\cov(\cM)}\ar[uu]\ar[ur]\ar[rrr]
&
&
& \sr{\sone(\B,\B)}{\cov(\cM)}\ar[uuuuu]\ar[ur]
}$
}
\caption{The extended Scheepers Diagram}\label{extSch}
\end{center}
\end{changemargin}
\end{figure}

Each of these properties, where
$\scrA,\scrB$ range over $\cO,\Lambda,\Omega,\Gamma$
or over $\cB$, $\BL$, $\BO$, $\BG$,
is either void or equivalent to one in Figure \ref{extSch}
(where an arrow denotes implication).
For these properties, $\cO$ can be replaced anywhere by
$\Lambda$ and $\cB$ by $\BL$ without changing the property \cite{coc1, coc2, CBC}.

The \emph{critical cardinality} of a property $\cI$
(relative to a space $X$) is
$$\non_X(\cI)=\min\{|Y| : Y\sbst X\mbox{ and }Y\nin\cI(X)\}.$$
The \emph{covering number} of $\cI$ (relative to $X$) is
$$\cov_X(\cI)=\min\{|\cF| : \cF\sbst \cI(X)\mbox{ and }\cup\cF= X\}.$$
Again, since we can work in $\bbR\sm\bbQ$, we remove the subscript $X$
from both notations. Below each property in Figure \ref{extSch} appears its
critical cardinality (these cardinals are well studied,
see \cite{BlassHBK}. By $\cM$ we always denote the ideal of meager, i.e.\ first category,
sets of real numbers).

$\sfin(\cO,\cO)$, $\ufin(\cO,\Gamma)$, $\sone(\cO,\cO)$
are the classical properties of
Menger, Hurewicz, and Rothberger (traditionally called $C'')$, respectively.
$\sone(\Omega,\Gamma)$ is the Gerlits-Nagy $\gamma$-property.
Additional properties in the diagram were studied by Arkhangel'ski\v{i},
Sakai, and others. Some of the properties are relatively new.

We also consider the following type of properties.
\bi
\itm[$\Split(\scrA,\scrB)$:] Every cover $\cU\in\scrA$ can be split
into two disjoint subcovers $\cV$ and $\cW$, each containing some element of $\scrB$ as a subset.
\ei
Here too, letting $\scrA,\scrB$ range over $\Lambda$, $\Omega$, $\Gamma$ or
$\BL$, $\BO$, $\BG$,
we get that some of the properties are trivial and several
equivalences hold among the remaining ones.
The surviving properties apper in the following diagram (where again the critical
cardinality appears below each property).
$${\scriptsize
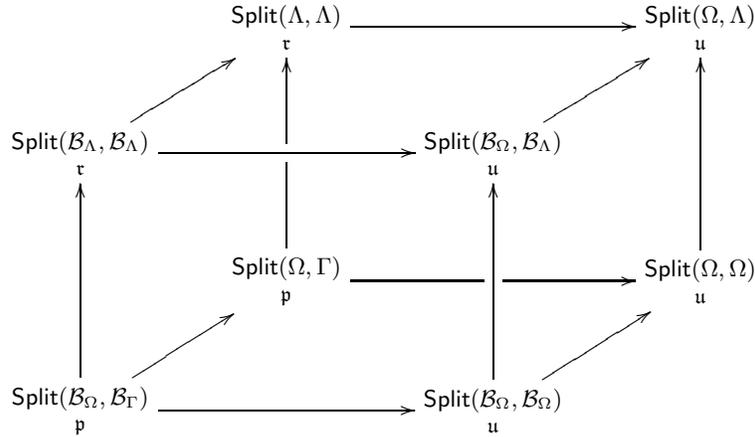
\xymatrix{
& \sr{\Split(\Lambda, \Lambda)}{\fr} \ar[rr] &  & \sr{\Split(\Omega, \Lambda)}{\fu}\\
\sr{\Split(\BL,\BL)}{\fr}\ar[ur]\ar[rr] & & \sr{\Split(\BO,\BL)}{\fu}\ar[ur]\\
& \sr{\Split(\Omega, \Gamma)}{\fp} \ar'[u][uu]\ar'[r][rr] & & \sr{\Split(\Omega, \Omega)}{\fu}\ar[uu]\\
\sr{\Split(\BO,\BG)}{\fp}\ar[rr]\ar[uu]\ar[ur] & & \sr{\Split(\BO,\BO)}{\fu}\ar[uu]\ar[ur]
}
}$$
No implication can be added to this diagram \cite{split}.
There are connections between the first and the second diagram,
e.g., $\Split(\Omega, \Gamma)=\sone(\Omega,\Gamma)$ \cite{split}, and
both $\ufin(\cO,\Gamma)$ and $\sone(\cO,\cO)$ imply $\Split(\Lambda,\Lambda)$.
Similarly, $\sone(\Omega,\Omega)$ implies $\Split(\Omega,\Omega)$ \cite{coc1}.
Similar assertions hold in the Borel case \cite{split}.

The situation becomes even more interesting when $\tau$-covers are incorporated into the framework.
We will introduce this notion later.

\section{Positive results}\label{positive}

\subsection{On the Scheepers diagram}

\begin{prop}[folklore]\label{easyadd}
Each property of the form $\Pi(\scrA,\cO)$ (or $\Pi(\scrA,\B)$),
$\Pi\in\{\sone,\sfin,\ufin\}$, is $\sigma$-additive.
\end{prop}
\begin{proof}
Let $A_1,A_2,\ldots$ be a partition of $\N$ into disjoint infinite
sets. Assume that $X_1,X_2\ldots$ satisfy $\Pi(\scrA,\cO)$. Assume
that $\cU_1,\cU_2,\ldots\in\scrA$ for $X=\Union_{k\in\N}X_k$.
For each $k$, use this property of $X_k$ to extract from the
sequence $\{\cU_n\}_{n\in A_k}$ the appropriate cover $\cV_k$ of
$X_k$. Then $\Union_{k\in\N}\cV_k$ is the desired cover of $X$.

The proof for $\Pi(\scrA,\B)$ is identical.
\end{proof}

\begin{prop}\label{imagele}
If $\I$ and $\cJ$ are collections of sets of reals such that:
\begin{quote}
$X\in\I$ if, and only if, for each Borel function $\Psi:X\to\R\sm\Q$
$\Psi[X]\in\cJ$.
\end{quote}
Then $\add(\cJ)\le\add(\I)$.
\end{prop}
\begin{proof}
Assume that $X_\alpha$, $\alpha<\kappa$, are members of $\I$ such that $X=\Union_{\alpha<\kappa}X_\alpha\nin\I$.
Take a Borel function $\Psi:X\to\R\sm\Q$ such that $\Psi[X]\nin\cJ$.
Then $\Psi[X]=\Union_{\alpha<\kappa}\Psi[X_\alpha]$.
\end{proof}

It is easy to see that for all $x,y\in\{\Gamma,\Omega,\cO\}$,
$X$ satisfies $\Pi(\B_x,\B_y)$ if, and only if,
every Borel image of $X$ satisfies $\Pi(x,y)$
(here $\B_\cO:=\B$) \cite{CBC, huremen}.
Using this and the facts that for each property $\cI$, $\add(\cI)$
is a regular cardinal satisfying $\add(\cI)\le\op{cf}(\non(\cI))$
and $\add(\cI)\le\cov(\cI)$, we have the following.

\begin{cor}
\mbox{}
\be
\itm $\add(\sone(\cO,\cO))\le\add(\sone(\B,\B))\le\cf(\cov(\cM))$,
\itm $\max\{\add(\sone(\Gamma,\Gamma)),\add(\ufin(\cO,\Gamma))\}\le
\add(\sone(\BG,\BG))\le\fb$,
\itm $\max\{\add(\sone(\Gamma,\cO)),\add(\sfin(\cO,\cO))\}\le
\add(\sfin(\B,\B))\le\allowbreak \cf(\fd)$,
\itm $\add(\sone(\Omega,\Gamma))\le\add(\sone(\BO,\BG))\le\fp$,
\itm $\max\{\add(\sone(\Gamma,\Omega)), \add(\sfin(\Gamma,\Omega)),
\add(\ufin(\cO,\Omega))\}\le$\\
$\le\add(\sone(\BG,\BO))\le\cf(\fd)$.\hfill\qed
\ee
\end{cor}

We now look for lower bounds on the additivity numbers.
Define a partial order $\le^*$ on $\NN$ by:
$$f\le^* g\quad \mbox{if}\quad f(n)\le g(n)\mbox{ for all but finitely many }n.$$
A subset of $\NN$ is called \emph{bounded} if it is bounded with respect to $\le^*$.
A subset $D$ of $\NN$ is \emph{dominating} if
for each $g\in\NN$ there exists $f\in D$ such that $g\le^* f$.

View $\N$ as a discrete topological space. The \emph{Baire space}
is the product space $\NN$.
Hurewicz (\cite{HURE27}, see also Rec\l{}aw \cite{RECLAW})
proved that a set of reals $X$
satisfies $\sfin(\cO,\cO)$ if, and only if,
every continuous image of $X$ in $\NN$ is not dominating.
Likewise, he showed that $X$ satisfies $\ufin(\cO,\Gamma)$ if,
and only if,
every continuous image of $X$ in $\NN$ is bounded.
Replacing ``continuous image'' by ``Borel image'' we get
characterizations of $\sfin(\B,\B)$ and $\sone(\BG,\BG)$,
respectively \cite{CBC}.
It is easy to see that
a union of less than $\fb$ many bounded subsets of $\NN$
is bounded, and
a union of less than $\fb$ many subsets of $\NN$ which are
not dominating is not dominating.

\begin{cor}\label{adds}
\mbox{}
\be
\itm $\add(\ufin(\cO,\Gamma))=\add(\sone(\BG,\BG))=\fb$;
\itm $\fb\le\add(\sfin(\cO,\cO))\le\add(\sfin(\B,\B))\le\cf(\fd)$. \hfill\qed
\ee
\end{cor}

Consider an unbounded subset $B$ of $\NN$ such that $|B|=\fb$, and define,
for each $f\in B$, $Y_f=\{g\in\NN : f\not\le^* g\}$. Then the
sets $Y_f$ are not dominating, but
$\Union_{f\in B} Y_f=\NN$: For each $g\in\NN$ there exists $f\in B$
such that $f\not\le^* g$, that is, $g\in Y_f$.
Thus the second assertion in Corollary \ref{adds} cannot be strengthened
in a trivial manner. We must work harder for that.

Let $\roth$ denote the collection of all infinite sets of natural numbers.
For $a,b\in\roth$, $a$ is an \emph{almost subset} of $b$, $a\as b$, if
$a\sm b$ is finite. A family $G\sbst\roth$ is \emph{groupwise dense} if
it contains all almost subsets of its elements, and for each
partition of $\N$ into finite intervals (i.e., sets of the form
$[m,k)=\{m,m+1,\ldots,k-1\}$), there is an infinite set of intervals
in this partition whose union is a member of $G$.

$\roth$ is a topological subspace of $\PN$, where the topology on
$\PN$ is defined by identifying it with the Cantor space $\Cantor$.
For each finite $F\sbst\N$ and each $n\in\N$, define
$$O_{F,n}=\{a\in\PN : a\cap[0,n)=F\}.$$
The sets $O_{F,n}$ form a clopen basis for the topology on $\PN$.

For $a\in\roth$, define an element $\Next{a}$ of $\NN$
by
$$\Next{a}(n) = \min\{k\in a : n<k\}$$
for each $n$.

\begin{thm}[{Tsaban-Zdomskyy \cite{MGD}}]\label{gd}
Assune that $X$ satisfies $\sfin\allowbreak (\cO,\cO)$.
Then for each continuous image $Y$ of $X$ in $\NN$,
the family
$$G=\{a\in\roth : (\forall f\in Y)\ \Next{a}\not\le^* f\}$$
is groupwise dense.
\end{thm}
\begin{proof}
Assume that $Y$ is a continuous image of $X$ in $\NN$.
Then $Y$ satisfies $\sfin(\cO,\cO)$.

\begin{lem}[folklore]\label{Ksigma}
Assume that $X$ satisfies $\sfin(\cO,\cO)$ and $K$ is $\sigma$-compact.
Then $X\x K$ satisfies $\sfin(\cO,\cO)$.
\end{lem}
\begin{proof}
This proof is as in  \cite{KMlinear}.
As $\sfin(\cO,\cO)$ is $\sigma$-additive, we may
assume that $K$ is compact.
Assume that $\cU_1,\cU_2,\dots$, are countable open covers of $X\x K$.
For each $n$, enumerate $\cU_n=\{U^n_m : m\in\N\}$.
For each $n$ and $m$ set
$$V^n_m = \left\{x\in X : \{x\}\x K\sbst \Union_{k\le m}U^n_k\right\}.$$
Then $\cV_n = \{V^n_m : m\in\N\}$ is an open cover of $X$.
As $X$ satisfies $\sfin(\cO,\cO)$, we can choose for each $n$ an $m_n$
such that $X=\Union_n\Union_{k\le m_n}V^n_k$.
By the definition of the sets $V^n_k$,
$X\x K\sbst\Union_n\Union_{k\le m_n}U^n_k$.
\end{proof}

By Lemma \ref{Ksigma}, $\PN\x Y$ satisfies $\sfin(\cO,\cO)$.

\begin{lem}[\cite{MGD}]\label{Clemma}
The set
$$C = \{(a,f)\in\roth\x\NN : \Next{a}\le^* f\}$$
is an $F_\sigma$ subset of $\PN\x\NN$.
\end{lem}
\begin{proof}
Note that
$$C = \Union_{m\in\N}\bigcap_{n\ge m}\{(a,f)\in \PN\x\NN : (n,f(n)]\cap a\neq\emptyset\}.$$
(The nonempty intersection for infinitely many $n$ allows the replacement of $\roth$ by $\PN$.)

For fixed $m$ and $n$, the set $\{(a,f)\in P(\N)\x\NN : (n,f(n)]\cap a\neq\emptyset\}$ is clopen:
Indeed, if $\lim_k(a_k,f_k)=(a,f)$ then for all large enough $k$, $f_k(n)=f(n)$,
and therefore for all larger enough $k$, $(n,f_k(n)]\cap a_k=(n,f(n)]\cap a$.
Thus, $(a_k,f_k)$ is in the set if, and only if, $(a,f)$ is in the set.
\end{proof}

As $\sfin(\cO,\cO)$ is $\sigma$-additive and hereditary for closed subsets,
we have by Lemma \ref{Clemma} that
$C\cap (\PN\x Y)$ satisfies $\sfin(\cO,\cO)$, and
therefore so does its projection $Z$ on the first coordinate.
By the definition of $Z$, $G=Z\comp$, the complement of $Z$ in $\roth$.
Note that $G$ contains all almost subsets of its elements.

For $a\in\roth$ and an increasing $h\in\NN$, define
$$a/h = \{n : a\cap \intvl{h}{n}\neq\emptyset\}.$$
For $S\sbst\roth$, define $S/h = \{a/h : a\in S\}$.

\begin{lem}[\cite{MGD}]\label{tala}
Assume that $G\sbst\roth$ contains all almost subsets of its elements.
Then: $G$ is groupwise dense if, and only if, for each increasing $h\in\NN$,
$G\comp/h\neq\roth$.
\end{lem}
\begin{proof}
For each increasing $h\in\NN$ and each $a\in\roth$,
$$\Union_{n\in a}\intvl{h}{n}\nin G\Iff \Union_{n\in a}\intvl{h}{n}\in G\comp
\Iff a\in G\comp/h.$$
The lemma follows directly from that.
\end{proof}

Assume that $G$ is not groupwise dense.
By Lemma \ref{tala}, there is an increasing $h\in\NN$ such that
$Z/h = G\comp/h=\roth$.
The natural mapping $\Psi:Z\to Z/h$ defined by $\Psi(a)=a/h$ is
a continuous surjection. It follows that $\roth$ satisfies $\sfin(\cO,\cO)$.
But this is absurd:
The image of $\roth$ in $\NN$, under the
continuous mapping assigning to each $a\in\roth$ its increasing enumeration,
is a dominating subset of $\NN$.
Thus, $\roth$ does not satisfy $\sfin(\cO,\cO)$ -- a contradiction.
\end{proof}

We obtain the promised improvement of Corollary \ref{adds}(2).

\begin{cor}[{Zdomskyy \cite{SF1, MGD}}]\label{addMen}
$\max\{\fb,\fg\}\le\add(\sfin(\cO,\cO))\le\add\allowbreak(\sfin(\B,\B))\le\cf(\fd)$.
\end{cor}
\begin{proof}
By Corollary \ref{adds}, we need only show that $\fg\le\add(\sfin(\cO,\cO))$.

Assume that $\kappa<\fg$ and for each $\alpha<\kappa$, $X_\alpha$ satisfies $\sfin(\cO,\cO)$,
and that $X=\Union_{\alpha<\kappa}X_\alpha$.
By the Hurewicz Theorem, it suffices to show that no continuous image of $X$ in
$\NN$ is dominating.
Indeed, assume that $\Psi:X\to\NN$ is continuous.
By Theorem \ref{gd}, for each $\alpha$ the family
$$G_\alpha =\{a\in\roth : (\forall f\in \Psi[X_\alpha])\ \Next{a}\not\le^* f\}$$
is groupwise dense.
Thus, there exists $a\in\bigcap_{\alpha<\kappa}G_\alpha$.
Then $\Next{a}$ witnesses that $\Psi[X]$ is not dominating.
\end{proof}

\begin{prob}
Is it consistent that $\max\{\fb,\fg\}<\add(\sfin(\cO,\cO))$?
\end{prob}

The methods used to obtain the last lower bound are similar to earlier methods of
Scheepers used to bound $\add(\sone(\Gamma,\Gamma))$ from below.
A family $\cD\sbst\roth$ is \emph{open} if it is closed under almost subsets.
It is \emph{dense} if for each $a\in\roth$ there is $d\in D$ such that
$d\as a$.
The \emph{density number} $\fh$ is the minimal cardinality of a collection of
open dense families in $\roth$ whose intersection is empty.
Identify $\roth$ with the increasing elements of $\NN$ by taking increasing enumerations.

\begin{thm}[Scheepers \cite{wqn}]\label{s1gg}
Assume that $X$ satisfies $\sone(\Gamma,\Gamma)$,
and $\cU_1,\cU_2,\ldots$ are open $\gamma$-covers of $X$.
For each $n$, enumerate $\cU_n = \{U^n_m : m\in\N\}$.
Then the family of all $a\in\roth$ such that
$\setseq{U^n_{a(n)}}$ is a $\gamma$-cover of $X$ is open dense.
\end{thm}
\begin{proof}
By standard arguments, we may assume that the given $\gamma$-covers are pairwise disjoint
(use the fact that any countable sequence of infinite sets can be refined to a countable
sequence of pairwise disjoint infinite sets.)

For each $n$ and $m$, define
$$V^n_m = U^1_m\cap U^2_m\cap\dots\cap U^n_m.$$
Fix any $a\in\roth$.
For each $n$, define
$$\cV_n = \{V^n_{a(m)} : m\ge n\}.$$
Then $\cV_n\in\Gamma$.
By $\sone(\Gamma,\Gamma)$, there is $f\in\NN$ such that $f(n)\ge n$ for all $n$,
and $\setseq{V^n_{a(f(n))}}\in\Gamma$.

Let $\tilde f$ be such that $\tilde f(1)=f(1)$, and for each $n\ge 1$, $\tilde f(n+1)=f(k)$
for some $k>n$ with $\tilde f(n)<f(k)$.
By the definition of the sets $V^n_m$,
$\setseq{V^n_{a(\tilde f(n))}}\in\Gamma$ as well.
Let $d\in\roth$ be such that $d(n)=a(\tilde f(n))$ for all $n$.
Then $d\sbst a$, and as $\setseq{V^n_{d(n)}}\in\Gamma$, we have
again by the definition of the sets $V^n_m$, that
$\setseq{V^n_{b(n)}}\in\Gamma$ for all $b\sbst d$.
In particular, $\setseq{U^n_{b(n)}}\in\Gamma$ for all $b\sbst d$.
\end{proof}

\begin{cor}[{Scheepers \cite{wqn}}]
$\fh\le\add(\sone(\Gamma,\Gamma))\le\add(\sone(\BG,\BG))\allowbreak\le\fb$.
\end{cor}
\begin{proof}
Fix $\kappa<\fh$ and assume that $X_\alpha$, $\alpha<\kappa$,
all satisfy $\sone(\Gamma,\Gamma)$.
Let $X=\Union_{\alpha<\kappa}X_\alpha$, and assume that for each $n$,
$\cU_n = \{U^n_m : m\in\N\}$ is an open $\gamma$-cover of $X$.

By Theorem \ref{s1gg}, for each $\alpha$ the family
$$\cD_\alpha = \{a\in\roth : \setseq{U^n_{a(n)}}\mbox{ is a $\gamma$-cover of }X\}$$
is open dense.
Take $a\in\bigcap_{\alpha<\kappa}\cD_\alpha$.
Then $\setseq{U^n_{a(n)}}$ is a $\gamma$-cover of $X$.
\end{proof}

\begin{prob}
Is it consistent that $\fh<\add(\sone(\Gamma,\Gamma))$?\footnote{Added after publication:
The answer is positive. Dow proved in \cite{Dow90} a theorem implying \cite{hH} that
$\sone(\Gamma,\Gamma)=\sone(\BG,\BG)$ (indeed, $\sone(\BG,\BG)=[\R]^{<\fb}$ in that model  \cite{BBC}).
In Laver's model, $\fh<\fb$ \cite{BlassHBK}. Apply Corollary \ref{adds}(1).}
\end{prob}

\begin{prob}
Is it consistent that $\add(\sone(\Gamma,\Gamma))<\fb$?
\end{prob}

We conclude the section with the following beautiful result.
Let $\cN$ denote the collection of Lebesgue null sets of reals.

\begin{thm}[{Carlson \cite{covM2}}]
$\add(\cN)\le\add(\sone(\cO,\cO))\le\add(\sone(\B,\B))\allowbreak\le\cf(\cov(\cM))$.
\end{thm}
\begin{proof}
The new ingredient is the first inequality.

\begin{lem}[Bartoszy\'nski \cite{barju}]\label{addNchar}
$\add(\cN)$ is the smallest cardinality of a family $F \sbst\NN$
such that there is no function $S :\N\to [\N]^{<\aleph_0}$ with $|S(n)| \leq n$ for all $n$,
such that $(\forall f \in F)(\forall^\infty n)\ f(n) \in S(n)$.\hfill\qed
\end{lem}

Assume that $\kappa<\add(\cN)$ and $X_\alpha$, $\alpha < \kappa$, satisfy $\sone(\cO,\cO)$.
Let $X = \bigcup_{\alpha < \kappa}X_\alpha$.
Assume that $\cU_n=\{U^n_m : m\in\N\}$, $n\in\N$, are open covers of $X$.
Let $r_n = 1+2+ \cdots +(n-1)$. For each $n$, let
$$\tilde\cU_n = \{\tilde U^n_s : s:[r_n,r_{n+1})\to\N\},$$
where $\tilde U^n_s = \bigcap_{k=r_n}^{r_{n+1}} U^k_{s(k)}.$
$\tilde\cU_n$ is an open cover of $X$.
For each $\alpha<\kappa$, as $X_\alpha$ satisfies $\sone(\cO,\cO)$,
there is $f_\alpha:\N\to\N^{<\alephes}$ such that $f(n)\in\N^n$ for each $n$,
and $\setseq{\tilde U^n_{f_\alpha(n)}}$ is a cover of $X_\alpha$.
By Lemma \ref{addNchar}, there is
$S :\N\to [\N^{<\aleph_0}]^{<\aleph_0}$ with $S(n)\in\N^n$ and $|S(n)| \leq n$ for all $n$,
such that
$$(\forall \alpha<\kappa)(\forall^\infty n)\ f_\alpha(n) \in S(n).$$
For each $n$, $S(n)$ contains at most $n$
sequences of length $n$. Let $g\in\NN$ be a function which agrees at least
once on the $n$-element interval $[r_n,r_{n+1})$ with each
of these sequences.
Then $\setseq{U^n_{g(n)}}$ is a cover of $X$.
\end{proof}

\subsection{On splitting properties}

\begin{thm}[\cite{split}]\label{add3}
$\Split(\BO,\BL)$ and $\Split(\Omega,\Lambda)$ are $\sigma$-additive.
\end{thm}
\begin{proof}
We will prove the open case. The Borel case is similar.
\begin{lem}[\cite{split}]
Assume that $\cU$ is a countable open $\omega$-cover of $Y$ and
that $X\subseteq Y$ satisfies $\Split(\Omega,\Lambda)$.
Then $\cU$ can be partitioned into two pieces $\cV$ and $\cW$
such that that $\cW$ is an $\omega$-cover of $Y$ and
$\cV$ is a large cover of $X$.
\end{lem}
\begin{proof}
First assume that there does not exist $U\in \cU$ with
$X\subseteq U$.  Then $\cU$ in an $\omega$-cover
of $X$. By the splitting property we can divide it into two
pieces each a large cover of $X$.  Since $\cU$ is an
$\omega$-cover of $Y$, one of the pieces is an $\omega$-cover of
$Y$, and the lemma is proved.
If there are only finitely many $U\in \cU$ with $X\subseteq U$,
then $\tilde \cU=\cU\sm\{U\in \cU : X\subseteq U\}$ is still
an $\omega$-cover of $Y$ and we can apply to it the above argument.

Thus, assume that there are infinitely many $U\in \cU$ with $X\subseteq U$.
Then take a partition of $\cU$ into two pieces such that each piece
contains infinitely many sets $U$ with $X\sbst U$.
One of the pieces must be an $\omega$-cover of $Y$.
\end{proof}
Assume that $Y=\Union_{n\in\N}X_n$ where each $X_n$ satisfies $\Split(\Omega,\Lambda)$,
and let $\cU_0$ be an open $\omega$-cover of $Y$.
Given $\cU_n$ an open $\omega$-cover of $Y$, apply the lemma twice to get a partition
$\cU_n=\cV_n^0\cup\cV_n^1\cup\cU_{n+1}$ such that
$\cU_{n+1}$ is an open $\omega$-cover of $Y$ and for each $i=0,1$, each element of
$X_n$ is contained in infinitely many $V\in \cV_n^i$.
Then the families $\cV^i=\Union_{n\in\N}\cV^i_n$, $i=0,1$, are disjoint large covers of $Y$
which are subcovers of $\cU_0$.
\end{proof}

Proposition \ref{imagele} implies the following.
\begin{cor} \mbox{}
\be
\itm $\add(\Split(\Lambda,\Lambda))\le\add(\Split(\BL,\BL))\le\cf(\fr)$,
\itm $\add(\Split(\Omega,\Lambda))\le\add(\Split(\BO,\BL))\le\cf(\fu)$,
\itm $\add(\Split(\Omega,\Omega))\le\add(\Split(\BO,\BO))\le\cf(\fu).$\hfill\qed
\ee
\end{cor}

However, $\Split(\Omega,\Omega)$ and $\Split(\BO,\BO)$ are not provably additive,
as we shall see in Section \ref{ConNegative}.

Concerning $\sigma$-additivity (or even just additivity, i.e.\ $\aleph_0$-additivity),
exactly one question remains open.
\begin{prob}\label{SpLamLamAdd}
Is $\Split(\Lambda,\Lambda)$ provably additive?
What about the Borel case?
\end{prob}

\section{Consistently negative results}\label{ConNegative}

Showing that a certain class is
not additive is apparently harder: All known results require axioms beyond ZFC.
This is often necessary, as will be seen in Section \ref{ConPositive}.

\subsection{On the Scheepers diagram}

For a sequence $\seq{X_n}$ of subsets of $X$, define $\liminf X_n =
\Union_m\bigcap_{n\ge m} X_n$. For a family $\cU$ of subsets of
$X$, $L(\cU)$ denotes its closure under the operation $\liminf$. A
set of reals $X$ has the property $(\delta)$ if for each open
$\omega$-cover $\cU$ of $X$, $X\in L(\cU)$.
The property $(\delta)$ was introduced by Gerlits and Nagy in \cite{GN},
where they showed that $\sone(\Omega,\Gamma)$ implies $(\delta)$.
The converse implication is still open.
It seems that the fact that $(\delta)$ is not provably additive
was not noticed before, but if follows from a combination of
results from \cite{JORG}, \cite{GM}, as we now show.

\begin{thm}\label{gammanotadd}
Assume \CH{}. Th\-en no class between $\sone(\BO,\BG)$ and $\sone(\Omega,\Gamma)$ or even $(\delta)$
(inclusive) is additive.
\end{thm}
\begin{proof}
By a theorem of Brendle \cite{JORG},
assuming CH there exists a set of reals $X$ of size continuum
such that all subsets of $X$ satisfy $\sone(\BO,\BG)$.

As $\sone(\BO,\BG)$ is closed under taking Borel (continuous is
enough) images, we may assume that $X\sbst (0,1)$. For
$Y\sbst (0,1)$, write $Y+1 = \{y+1 : y\in Y\}$ for the translation of $Y$ by
$1$.
The following is essentially proved in Theorem 5 of Galvin and Miller's paper \cite{GM}.
\begin{lem}
If $Y\sbst X\sbst (0,1)$ and $Z=(X\sm Y)\cup (Y+1)$ has property $(\delta)$,
then $Y$ is a Borel subset of $X$.
\end{lem}
\begin{proof}
Let
$$\cU=\{U\cup (V+1) : \mbox{open }U,V\sbst (0,1),\ \cl{U}\cap \cl{V}=\emptyset\}.$$
$\cU$ is an open $\omega$-cover of $Z$.
If $U_n\cap V_n=\emptyset$ for all $n$, then the sets
$U=\Union_m\bigcap_{n\ge m} U_n$ and $V=\Union_m\bigcap_{n\ge m}V_n$ are disjoint,
and $\Union_m\bigcap_{n\ge m} U_n\cup (V_n+1)=U\cup (V+1)$.
It follows by transfinite induction, each element in $L(\cU)$ has the form
$U\cup (V+1)$ where $U,V$ are disjoint Borel subsets of $Z$.
Thus, if $Z\in L(\cU)$, there are such $U$ and $V$ with $Z=U\cup (V+1)$.
It follows that $Y=V\cap X$ is a Borel subset of $X$.
\end{proof}
As $|X|=\fc$ and only $\fc$ many out of the $2^\fc$ many
subsets of $X$ are Borel, there exists a subset $Y$ of $X$ which
is not Borel.
It follows that $(X\sm Y)\cup (Y+1)$ does not have the property $(\delta)$
(and, in particular, does not have the property $\sone(\Omega,\Gamma)$).
But by the choice of $X$, both $X\sm Y$ and $Y$ (and therefore also $Y+1$)
satisfy $\sone(\BO,\BG)$.
\end{proof}

Except for the $(\delta)$ part,
Theorem \ref{gammanotadd} was proved in \cite{tautau}.
The extension to $(\delta)$ was noticed by Miller (personal communication).

We next show that if $\cov(\cM)=\fc$ (in particular, assuming \CH{}), then
no class between $\sone(\BO,\BO)$ and $\ufin(\cO,\Omega)$ (inclusive) is additive.

For clarity of exposition, we will first treat the
open case, and then explain how to modify the constructions in order
to cover the Borel case.

For convenience, we will work in $\NZ$ (with pointwise addition), which is
homeomorphic to $\R\sm\Q$.
The notions that we will use are topological,
thus the following constructions can be translated to constructions
in $\R\sm\Q$.

A collection $\cJ$ of sets of reals is \emph{translation invariant}
if for each real $x$ and each $X\in\cJ$, $x+X\in\cJ$.
$\cJ$ is \emph{negation invariant} if for each $X\in\cJ$, $-X\in\cJ$
as well.
For example,
$\cM$ and $\Null$ are negation and translation invariant (and there are many
more examples).
\begin{lem}[folklore]\label{x+y=z}
If $\cJ$ is negation and translation invariant and
if $X$ is a union of less than $\cov(\cJ)$ many
elements of $\cJ$,
then for each $x\in\NZ$ there exist $y,z\in\NZ\sm X$ such that $y+z=x$.
\end{lem}
\begin{proof}
$(x - X) \cup X$ is a union of less than
$\cov(\cJ)$ many elements of $\cJ$.
Thus we can choose an element
$y\in\NZ\sm ((x-X)\cup X)=(x - \NZ\sm X)\cap (\NZ\sm X)$; therefore
there exists $z\in \NZ\sm X$
such that $x-z = y$, that is, $x = y+z$.
\end{proof}

A set of reals $L$ is \emph{$\kappa$-Luzin} if $|L|\ge\kappa$ and for
each meager set $M$, $|L\cap M|<\kappa$.

The following result was obtained independently by
many authors: A comment on the top of Page 205 of \cite{coc2} (without proof);
Theorem 13 of \cite{lengthdiags} (under \CH);
Section 3 of \cite{KMlinear}; Theorem 4 of \cite{huremen2};
Theorem 2 of \cite{gamma7} (under \CH).

\begin{prop}[folklore]\label{opennotadd}
Assume that $\cov(\cM)=\fc$.
Then there exist
$\fc$-Luzin subsets $L_0$ and $L_1$ of $\NZ$
satisfying $\sone(\Omega,\Omega)$,
such that $L_0+L_1=\NZ$.
\end{prop}
\begin{proof}
Assume that $\cov(\cM)=\fc$.
Let $\{y_{\alpha}:\alpha<\fc\}$ enumerate $\NZ$;
let $\{M_{\alpha}:\alpha<\fc\}$
enumerate all $F_\sigma$ meager sets in $\NZ$ (observe that this
family is cofinal in $\cM$),
and let $\{\seq{\cU^\alpha_n} : \alpha<\fc\}$
enumerate all countable sequences of countable families of open sets.

Fix a countable dense subset $Q\sbst\NZ$.
We construct $L_0=\{x^0_\beta : \beta<\fc\}\cup Q$ and
$L_1=\{x^1_\beta : \beta<\fc\}\cup Q$ by induction on
$\alpha<\fc$. During the construction, we make an inductive
hypothesis and verify that it remains true after making
the inductive step.

At stage $\alpha\ge 0$ set
\begin{eqnarray*}
X^0_\alpha & = & \{x^0_\beta : \beta<\alpha\}\cup Q\\
X^1_\alpha & = & \{x^1_\beta : \beta<\alpha\}\cup Q
\end{eqnarray*}
and consider the sequence $\seq{\cU^\alpha_n}$.
For each $i<2$, do the following.
Call $\alpha$ \emph{$i$-good} if
for each $n$ $\cU^\alpha_n$ is an $\omega$-cover of $X^i_\alpha$.
Assume that $\alpha$ is $i$-good.
Since $\cov(\cM)=\non(\sone(\Omega,\Omega))$ \cite{coc2}
and we assume that $\cov(\cM)=\fc$,
there exist elements
$U^{\alpha,i}_n\in\cU^\alpha_n$ such that $\seq{U^{\alpha,i}_n}$ is
an $\omega$-cover of $X^i_\alpha$.
We make the \emph{inductive hypothesis} that
for each $i$-good $\beta<\alpha$,
$\seq{U^{\beta,i}_n}$ is an $\w$-cover of $X^i_\alpha$.
For each finite $F\sbst X^i_\alpha$, and each $i$-good $\beta\le\alpha$,
define
$$G_i^{F,\beta}=\bigcup\{U^{\beta,i}_n : n\in\N,\ F\sbst U^{\beta,i}_n\}.$$
Then $Q\sbst G_i^{F,\beta}$ and thus $G_i^{F,\beta}$ is open and dense.

Set
$$Y_\alpha=\Union_{\beta<\alpha}M_\beta\cup
\Union\left\{\NZ\sm G_i^{F,\beta} : i<2, \beta\le\alpha\mbox{ $i$-good},F\sbst X^i_\alpha\mbox{ finite}\right\}.$$
Then $Y_\alpha$ is a union of less than $\cov(\cM)$ many meager sets, thus
by Lemma \ref{x+y=z}
we can pick $x^0_\alpha,x^1_\alpha\in\NZ\sm Y_\alpha$ such that $x^0_\alpha+x^1_\alpha=y_\alpha$.
To see that the inductive hypothesis is preserved, observe that
for each finite $F\sbst X^i_\alpha$ and $i$-good $\beta\le\alpha$,
$x^i_\alpha\in G_i^{F,\beta}$ and therefore $F\cup\{x^i_\alpha\}\sbst U^{\beta,i}_n$ for
some $n$.

Clearly $L_0$ and $L_1$ are $\fc$-Luzin sets, and $L_0+L_1=\NZ$.
It remains to show that $L_0$ and $L_1$ satisfy
$\sone(\Omega,\Omega)$.

Fix $i<2$. Consider, for each $\beta<\fc$, the
sequence $\seq{\cU^{\beta}_n}$.
If all members of that sequence are $\w$-covers of $L_i$,
then in particular they $\w$-cover $X^i_\beta$ (that is, $\beta$ is $i$-good).
By the inductive hypothesis, $\setseq{U^{\beta,i}_n}$
is an $\w$-cover of $X^i_\alpha$ for each $\alpha<\fc$, and therefore
an $\w$-cover of $L_i$.
\end{proof}

For a finite subset $F$ of $\NN$, define $\max(F)\in\NN$ to be the
function $g$ such that $g(n)=\max\{f(n) : f\in F\}$ for each $n$.
A subset $Y$ of $\NN$, is \emph{finitely-dominating} if the collection
$$\maxfin(Y):=\{\max(F) : F\mbox{ is a finite subset of }Y\}$$
is dominating.

\begin{thm}[{Tsaban \cite{huremen}}, Eisworth-Just \cite{gamma7}]\label{char1}
For a set of reals $X$, the following are equivalent:
\be
\itm $X$ satisfies $\ufin(\cO,\Omega)$;
\itm No continuous image of $X$ in $\NN$ is finitely-dominating.\hfill\qed
\ee
\end{thm}

A subset $Y$ of $\NN$ is \emph{$k$-dominating} if for each $g\in\NN$
there exists a $k$-element subset $F$ of $Y$ such that
$g\le^*\max(F)$ \cite{BlassNew}.
Clearly each $k$-dominating subset of $\NN$
is also finitely dominating.

Proposition \ref{opennotadd} and Theorem \ref{char1} imply
that no property between $\sone(\Omega,\Omega)$ and
$\ufin(\cO,\Omega)$ (inclusive) is provably additive.
Surprisingly, this was only observed in \cite{huremen2}.\footnote{Indeed,
in \cite{lengthdiags} Scheepers points out that Proposition \ref{opennotadd}
implies that no class between $\sone(\Omega,\Omega)$ and $\sfin(\Omega,\Omega)$
is provably additive. The missing ingredient to upgrade to $\ufin(\cO,\Omega)$ was
Theorem \ref{char1}.}

\begin{cor}[Bartoszy\'nski-Shelah-Tsaban \cite{huremen2}]
Assume that $\cov(\cM)\allowbreak =\fc$.
Then there exist $\fc$-Luzin subsets $L_0$ and $L_1$ of $\NZ$
satisfying $\sone(\Omega,\Omega)$, such that the $\fc$-Luzin set $L_0\cup L_1$ is
$2$-dominating. In particular, $L_0\cup L_1$ does not satisfy $\ufin(\cO,\Omega)$.
\end{cor}
\begin{proof}
Let $L_0,L_1$ be as in Proposition \ref{opennotadd}.
As $L_0+L_1=\NZ$ and in general $(f+g)/2\le \max\{f,g\}$ for all $f,g\in\NZ$,
we have that $L_0\cup L_1$ is $2$-dominating.
By Theorem \ref{char1}, the continuous image $\{|f| : f\in L_0\cup L_1\}$ of $L_0\cup L_1$
does not satisfy $\ufin(\cO,\Omega)$.
\end{proof}

We now treat the Borel case.
\begin{thm}[Bartoszy\'nski-Shelah-Tsaban \cite{huremen2}]\label{Borelnotadd}
Assume that $\cov(\cM)\allowbreak =\fc$.
Then there exist
$\fc$-Luzin subsets $L_1$ and $L_2$ of $\NZ$
satisfying $\sone(\BO,\BO)$,
such that for each $g\in\NZ$
there are $f_0\in L_0, f_1\in L_1$ satisfying $f_1(n)+f_2(n)=g(n)$
for all but finitely many $n$.

In particular, the $\fc$-Luzin set $L_0\cup L_1$
is $2$-dominating, and consequently does not satisfy $\ufin(\cO,\Omega)$.
\end{thm}
\begin{proof}
We follow the proof steps of Proposition \ref{opennotadd}.
The major problem is that here the sets $G_i^{F,\beta}$
need not be comeager. In order to overcome this,
we will consider only $\w$-covers where these sets are
guaranteed to be comeager, and make sure that it
is enough to restrict attention to this special sort of
$\w$-covers.
The following definition is essentially due to \cite{CBC},
but with a small twist that makes it work.
\begin{defn}[\cite{huremen2}]
A cover $\cU$ of $X$ is \emph{$\w$-fat} if for each finite $F\sbst X$
and each finite family $\cF$ of
nonempty open sets, there exists $U\in\cU$
such that $F\sbst U$ and for each $O\in\cF$, $U\cap O$ is not meager.
(Thus each $\w$-fat cover is an
$\w$-cover.) Let $\BOfat$ denote the collection of countable $\w$-fat
Borel covers of $X$.
\end{defn}

\begin{lem}[\cite{huremen2}]\label{fatlemma}
Assume that $\cU$ is a countable collection of Borel sets of reals.
Then $\cup\cU$ is comeager if, and only if,
for each nonempty basic open set $O$ there exists $U\in\cU$
such that $U\cap O$ is not meager.
\end{lem}
\begin{proof}
($\Rightarrow$) Assume that $O$ is a nonempty basic open set.
Then $\cup\cU\cap O=\Union\{U\cap O : U\in\cU\}$ is a countable union
which is not meager. Thus there exists $U\in\cU$
such that $U\cap O$ is not meager.

($\Leftarrow$) Set $B=\cup\cU$.
As $B$ is Borel, it has the Baire property.
Let $O$ be an open set and $M$ be a meager set such that
$B=(O\sm M)\cup (M\sm O)$.
For each basic open set $G$, $B\cap G$ is not meager, thus
$O\cap G$ is not meager as well.
Thus, $O$ is open dense. As $O\sm M\sbst B$, we have that
$\R\sm B\sbst(\R\sm O)\cup M$ is meager.
\end{proof}

\begin{cor}[\cite{huremen2}]\label{addelement}
Assume that $\cU$ is an $\w$-fat cover of some set $X$.
Then:
\be
\itm For each finite $F\sbst X$ and finite family $\cF$ of
nonempty \emph{basic} open sets, the set
$$\Union\left\{U\in\cU : F\sbst U\mbox{ and for each $O\in\cF$, }U\cap O\nin\cM\right\}$$
is comeager.
\itm For each element $x$ in the intersection of all sets of this form,
$\cU$ is an $\w$-fat cover of $X\cup\{x\}$.
\ee
\end{cor}
\begin{proof}
Write
$$\cV_{F,\cF} = \{U\in\cU : F\sbst U\mbox{ and for each $O\in\cF$, }U\cap O\nin\cM\}.$$
(1) Assume that $G$ is a nonempty open set. As $\cU$ is $\w$-fat and
the family $\cF\cup\{G\}$ is finite, there exists $U\in\cV_{F,\cF}$ such that
$U\cap G$ is not meager.
By Lemma \ref{fatlemma}, $\cup\cV_{F,\cF}$ is comeager.

(2) Assume that $F$ is a finite subset of $X\cup\{x\}$ and
$\cF$ is a finite family of nonempty basic open sets.
As $x\in\cup\cV_{F\sm\{x\},\cF}$, there exists $U\in\cU$ such that $x\in U$,
$F\sm\{x\}\sbst U$ (thus $F\sbst U$), and for each $O\in\cF$, $U\cap O$ is not meager.
\end{proof}

\begin{lem}[\cite{huremen2}]\label{s1fatfat}
If $|X|<\cov(\cM)$, then $X$ satisfies $\sone(\BOfat,\BOfat)$.
\end{lem}
\begin{proof}
Assume that $|X|<\cov(\cM)$, and let $\seq{\cU_n}$ be a sequence of
countable Borel $\w$-fat covers of $X$.
Enumerate each cover $\cU_n$ by $\{U^n_k\}_{k\in\N}$.
Let $\seq{A_n}$ be a partition of $\N$ into infinitely many
infinite sets. For each $m$, let $a_m\in\NN$ be an increasing
enumeration of $A_m$.
Let $\seq{\cF_n}$ be an enumeration of all finite families of
nonempty basic open sets.

For each finite subset $F$ of $X$ and each $m$ define a function
$\Psi^m_F\in\NN$ by
$$\Psi^m_F(n)=\min\{k : F\sbst U^{a_m(n)}_k\mbox{ and for each $O\in\cF_m$, }U^{a_m(n)}_k\cap O\nin\cM\}$$
Since there are less than $\cov(\cM)$ many functions $\Psi^m_F$,
there exists by \cite{covM} a function $f\in\NN$ such that
for each $m$ and $F$, $\Psi^m_F(n)=f(n)$ for infinitely many $n$.
Consequently, $\cV=\{U^{a_m(n)}_{f(n)} : m,n\in\N\}$ is an $\w$-fat cover of $X$.
\end{proof}

The following lemma justifies our focusing on $\w$-fat covers.

\begin{lem}[\cite{huremen2}]\label{denselusin}
Assume that $L$ is a set of reals such that for each nonempty basic
open set $O$, $L\cap O$ is not meager.
Then every countable Borel $\w$-cover $\cU$ of $L$ is an $\w$-fat cover of $L$.
\end{lem}
\begin{proof}
Assume that $\cU$ is a countable collection of Borel sets which is
not an $\w$-fat cover of $L$.
Then there exist a finite set $F\sbst L$ and nonempty open sets $O_1,\dots,O_k$
such that for each $U\in\cU$ containing $F$, $U\cap O_i$ is meager for some $i$.
For each $i=1,\dots,k$ let
$$M_i = \Union\left\{U\in\cU : F\sbst U\mbox{ and }U\cap O_i\in\cM\right\}.$$
Then $M_i\cap O_i$ is meager,
thus there exists $x_i\in (L\cap O_i)\sm M_i$.
Then $F\cup\{x_1,\dots,x_k\}$ is not covered by any $U\in\cU$.
\end{proof}

Let $\NZ=\{y_{\alpha}:\alpha<\fc\}$,
$\{M_{\alpha}:\alpha<\fc\}$ be all $F_\sigma$ meager subsets of $\NZ$, and
$\{\seq{\cU^\alpha_n} : \alpha<\fc\}$ be all sequences of countable families of
Borel sets.
Let $\{O_k : k\in\N\}$ and $\{\cF_m: m\in\N\}$ be all nonempty basic open sets
and all finite families of nonempty basic open sets, respectively, in $\NZ$.

We construct $L_i=\{x^i_\beta : \beta<\fc\}$, $i=0,1$,
by induction on $\alpha<\fc$ as follows.
At stage $\alpha\ge 0$ set
$X^i_\alpha  = \{x^i_\beta : \beta<\alpha\}$
and consider the sequence $\seq{\cU^\alpha_n}$.
Say that $\alpha$ is $i$-good if for each $n$
$\cU^\alpha_n$ is an $\w$-fat cover of $X^i_\alpha$.
In this case,
by Lemma \ref{s1fatfat} there exist elements
$U^{\alpha,i}_n\in\cU^\alpha_n$ such that $\seq{U^{\alpha,i}_n}$ is
an $\w$-fat cover of $X^i_\alpha$.
We make the inductive hypothesis that
for each $i$-good $\beta<\alpha$,
$\seq{U^{\beta,i}_n}$ is an $\w$-fat cover of $X^i_\alpha$.
For each finite $F\sbst X^i_\alpha$, $i$-good $\beta\le\alpha$,
and $m$ define
$$G_i^{F,\beta,m}=\Union\left\{U^{\beta,i}_n : F\sbst U^{\beta,i}_n
\mbox{ and for each $O\in\cF_m$, }U^{\beta,i}_n\cap O\nin\cM\right\}.$$
By Corollary \ref{addelement}(1), $G_i^{F,\beta,m}$ is comeager.
Set
$$Y_\alpha=\Union_{\beta<\alpha}M_\beta\cup
\Union\left\{\NZ\sm G_i^{F,\beta,m} :
\begin{matrix}
i<2,\beta\le\alpha\mbox{ $i$-good},\\
m\in\N,F\sbst X^i_\alpha\mbox{ Finite}
\end{matrix}
\right\}.$$
and $Y_\alpha^* = \{x\in\NZ : (\E y\in Y_\alpha)\ x=^* y\}$ (where $x =^* y$ means that $x(n)=y(n)$ for
all but finitely many $n$.) Then $Y_\alpha^*$ is a union of less
than $\cov(\cM)$ many meager sets. Use Lemma \ref{x+y=z} to pick
$x^0_\alpha,x^1_\alpha\in\NZ\sm Y_\alpha^*$ such that
$x^0_\alpha+x^1_\alpha=y_\alpha$. Let $k = \alpha \bmod \omega$,
and change a finite initial segment of $x^0_\alpha$ and
$x^1_\alpha$ so that they both become members of $O_k$. Then
$x^0_\alpha,x^1_\alpha\in O_k\sm Y_\alpha$, and
$x^0_\alpha+x^1_\alpha =^* y_\alpha$. By Corollary
\ref{addelement}(2), the inductive hypothesis is preserved.

Thus each $L_i$ satisfies $\sone(\BOfat,\BOfat)$
and its intersection with each nonempty basic open set
has size $\fc$.
By Lemma \ref{denselusin}, $\BOfat=\BO$ for $L_i$.
Finally, $L_0+L_1$ is dominating, so $L_0\cup L_1$ is $2$-dominating.
\end{proof}

Thus, no class between $\sone(\BO,\BO)$ and $\ufin(\cO,\Omega)$ (inclusive) is
provably additive.

\begin{rem}
As $\non(\ufin(\cO,\Omega))=\fd$, a natural question is whether
the method of Proposition \ref{opennotadd} can be generalized
to work for $\ufin(\cO,\Omega)$ under the weaker assumption
$\fd=\fc$.
By the forthcoming Theorem \ref{Sch_u<g}, such a trial is doomed to fail,
since $\fu<\fg$ implies that $\fg=\fd=\fc$.
\end{rem}

\subsection{On splitting properties}

It is well known that nonprincipal ultrafilters on $\N$
do not have the Baire property, and in particular are nonmeager \cite{barju}.
We can prove more than that.

\begin{lem}[Shelah \cite{split}]\label{UminusM}
Assume that $U$ is a nonprincipal ultrafilter on $\N$ and that
$M\sbst\roth$ is meager. Then $U\sm M$ is a subbase for $U$.
In fact, for each $a\in U$ there exist $a_0,a_1\in U\sm M$ such
that $a_0\cap a_1\sbst a$.
\end{lem}
\begin{proof}
Recall that $\roth$ is a subspace of $P(\N)$ whose topology is
defined by its identification with $\Cantor$.
It is well known \cite{barju, BlassHBK} that for each meager subset $M$ of $\Cantor$
there exist $x\in\Cantor$ and an increasing
$h\in\NN$ such that
$$M\sbst\{y\in\Cantor : (\forall^\oo n)\ y\upharpoonright \intvl{h}{n}\neq x\upharpoonright \intvl{h}{n}
\}.$$
(The set on the right hand side is also meager.)
Translating this to the language of $\roth$, we get that
for each $n$ there exist disjoint sets $I^n_0$ and $I^n_1$ satisfying
$I^n_0\cup I^n_1 = \intvl{h}{n}$,
such that
\begin{equation}\label{Mcof}
M\sbst\{y\in\roth : (\forall^\oo n)\ y\cap I^n_0\neq\emptyset\mbox{ or } I^n_1\not\sbst y\}.
\end{equation}

Assume that the sets $I^n_0,I^n_1$, $n\in\N$, are chosen as in
\eqref{Mcof}.
Let $a$ be an infinite co-infinite subset of $\N$.
Then either $x=\Union_{n\in a}\intvl{h}{n}\nin U$, or else
$x=\Union_{n\in \N\sm a}\intvl{h}{n}\nin U$.
We may assume that the former case holds.
Split $a$ into two disjoint infinite sets $a_1$ and $a_2$.
Then $x_i=\Union_{n\in a_i}\intvl{h}{n}\nin U$ ($i=0,1$).

Assume that $b\in U$. Then $\tilde b = b\sm x = b\cap(\N\sm x)\in U$.
Define sets $y_1,y_2\in U\sm M$ as follows.
\begin{eqnarray*}
y_1 & = & \tilde b \cup \Union_{n\in a_2}I^n_1\\
y_2 & = & \tilde b \cup \Union_{n\in a_1}I^n_1
\end{eqnarray*}
By (\ref{Mcof}), $y_1,y_2\nin M$.
As $y_1,y_2\spst\tilde b$, $y_1,y_2\in U$.
Now, $y_1\cap y_2=\tilde b\sbst b$.
\end{proof}

\begin{thm}[Tsaban \cite{split}]\label{notadd}
Assume that $\add(\cM)=\fc$. Then there exist two $\fc$-Luzin sets
$L_0$ and $L_1$ such that:
\be
\itm $L_0,L_1$ satisfy $\sone(\BO,\BO)$,
\itm $L=L_0\cup L_1$ satisfies $\Split(\BL,\BL)$; and
\itm $L=L_0\cup L_1$ does not satisfy $\Split(\Omega,\Omega)$.
\ee
\end{thm}
\begin{proof}
We follow the footsteps of the proof of Theorem \ref{Borelnotadd}.
Let $U=\{a_\alpha : \alpha<\fc\}$ be a nonprincipal ultrafilter on $\N$.
Let $\{M_{\alpha}:\alpha<\fc\}$ enumerate all $F_\sigma$ meager sets in $\roth$,
and $\{\seq{\cU^\alpha_n} : \alpha<\fc\}$
enumerate all countable sequences of countable families of Borel sets in $\roth$.
Let $\{O_i : i\in\N\}$ and $\{\cF_i : i\in\N\}$ enumerate all nonempty basic open sets
and finite families of nonempty basic open sets, respectively, in $\roth$.

We construct $L_i=\{a^i_\beta : \beta<\fc\}$, $i=0,1$,
by induction on $\alpha<\fc$ as follows.
At stage $\alpha\ge 0$ set
$X^i_\alpha  = \{a^i_\beta : \beta<\alpha\}$
and consider the sequence $\seq{\cU^\alpha_n}$.
Say that $\alpha$ is $i$-good if for each $n$
$\cU^\alpha_n$ is an $\w$-fat cover of $X^i_\alpha$.
In this case,
by the above remarks there exist elements
$U^{\alpha,i}_n\in\cU^\alpha_n$ such that $\seq{U^{\alpha,i}_n}$ is
an $\w$-fat cover of $X^i_\alpha$.
We make the inductive hypothesis that
for each $i$-good $\beta<\alpha$,
$\seq{U^{\beta,i}_n}$ is an $\w$-fat cover of $X^i_\alpha$.
For each finite $F\sbst X^i_\alpha$, $i$-good $\beta\le\alpha$,
and $m$ define
$$G_i^{F,\beta,m}=\Union\left\{U^{\beta,i}_n : F\sbst U^{\beta,i}_n\mbox{ and }
(\forall O\in\cF_m)\ U^{\beta,i}_n\cap O\nin\cM\right\}.$$
By the inductive hypothesis, $G_i^{F,\beta,m}$ is comeager.
Set
$$Y_\alpha=\Union_{\beta<\alpha}M_\beta\cup
\Union\left\{\roth\sm G_i^{F,\beta,m} :
\begin{matrix}
i<2,\beta\le\alpha\mbox{ $i$-good},\\
m\in\N,F\sbst X^i_\alpha\mbox{ Finite}
\end{matrix}
\right\},$$
and $Y_\alpha^* = \{x\in\roth : (\E y\in Y_\alpha)\ x=^* y\}.$
(Here $x =^* y$ means that $x\as y$ and $y\as x$.)
$Y_\alpha^*$ is a union of less than $\add(\cM)$ many meager sets,
and is therefore meager.
Use Lemma \ref{UminusM}
to pick $a^0_\alpha,a^1_\alpha\in U\sm Y_\alpha^*$ such that
$a^0_\alpha\cap a^1_\alpha \as a_\alpha$.
Let $k = \alpha \bmod \omega$, and change finitely many elements
of $a^0_\alpha$ and $a^1_\alpha$ so that
they both become members of $O_k$.
Then $a^0_\alpha,a^1_\alpha\in (U\cap O_k)\sm Y_\alpha$,
and $a^0_\alpha \cap a^1_\alpha \as a_\alpha$.
Observe that the inductive hypothesis remains true for $\alpha$.
This completes the construction.

Clearly $L_0$ and $L_1$ are $\fc$-Luzin sets
and $L_0\cup L_1$ is a subbase for $U$.
We made sure that for each nonempty basic open set $G$,
$|L_0\cap G|=|L_1\cap G|=\fc$,
thus $\BO=\BOfat$ for $L_0$ and $L_1$.
By the construction, $L_0,L_1\in\sone(\BOfat,\BOfat)$.

As we assume that $\add(\cM)=\fc$, every $\fc$-Luzin set (in particular, $L_0\cup L_1$) satisfies
$\sone(\B,\B)$ \cite{CBC}, and therefore also $\Split(\BL,\BL)$.

\begin{lem}[Just-Miller-Scheepers-Szeptycki \cite{coc2}]
If there is a continuous image of $X$ in $\roth$ that is a subbase for
a nonprincipal ultrafilter on $\N$, then $X$ does not satisfy $\Split(\Omega,\Omega)$.\hfill\qed
\end{lem}

As $L_0\cup L_1$ is a subbase for a nonprincipal ultrafilter on $\N$,
it does not satisfy $\Split(\Omega,\Omega)$.
\end{proof}

It follows that no property between $\sone(\BO,\BO)$ and $\Split(\Omega,\Omega)$ is
provably additive.

\section{Consistently positive results}\label{ConPositive}

\subsection{On the Scheepers diagram}

\begin{thm}[folklore]
It is consistent that all classes between $\sone(\Omega,\allowbreak\Gamma)$ and $\sone(\cO,\cO)$
(inclusive) are $\sigma$-additive.
\end{thm}
\begin{proof}
As $\sone(\cO,\cO)$ implies strong measure zero,
Borel's Conjecture (which asserts that every strong measure zero set is countable)
implies that all elements of $\sone(\cO,\cO)$ are countable, and thus
all classes below $\sone(\cO,\cO)$ are $\sigma$-additive.
Borel's Conjecture was proved consistent by Laver \cite{LAVER}.
\end{proof}

A variant of Borel's Conjecture for $\ufin(\cO,\Omega)$ is
false \cite{coc2, wqn, ideals, SFH}.
However, we have the following.

\begin{thm}[Bartoszy\'nski-Shelah-Tsaban \cite{huremen2}, Zsomskyy \cite{SF1, combimgs}]\label{Sch_u<g}
If $\fu<\fg$, then $\add(\ufin(\cO,\Omega))=\add(\sone(\BG,\BO))=\fc$.
\end{thm}
\begin{proof}
In \cite{SF1, combimgs} it is proved that
$\fu<\fg$ implies that $\ufin(\cO,\Omega)=\sfin(\cO,\cO)$, and
the same assertion holds in the Borel case.
The theorem follows from Corollary \ref{addMen}, together with the
fact that $\fu<\fg$ implies that $\fg=\fc$ \cite{BlassHBK}.
\end{proof}

In the remainder of this section we will show that
$\sigma$-additivity of $\ufin(\cO,\Omega)$ (and $\sone(\BG,\BO)$)
actually follow from the weaker axiom NCF, and that
a suitable combinatorial version of this assertion actually
gives a characterization of NCF.

In Theorem \ref{char1}, $\NN$ can be replaced by
$\NNup$ -- the \emph{(strictly) increasing} elements of $\NN$.
To see this, note that the function
$\Phi:\NN\to \NNup$
defined by
$$\Phi(f)(n)=n+f(0)+f(1)+\ldots+f(n)$$
is a homeomorphism which preserves finite-dominanace
in both directions.

We now consider the purely combinatorial counterpart of the
question whether $\ufin(\cO,\Omega)$ is additive.
Let $\Dfin$ denote the collection of subsets
of $\NNup$ which are not finitely-dominating.
By the previous comment,
$$\add(\Dfin)\le \add(\ufin(\cO,\Omega))\le \add(\sone(\BG,\BO)).$$

Recall that for an increasing $h\in\NN$ and a filter $\cF\sbst\roth$,
$$\cF/h= \{a/h : a\in\cF\} = \left\{ a : \Union_{n\in a}\intvl{h}{n}\in\cF\right\}.$$
(The first equality is the definition; the second an easy fact.)
If $\cF$ is an ultrafilter, then so is $\cF/h$.
We say that filters $\cF_1$ and $\cF_2$ on $\N$ are \emph{compatible in the Rudin-Keisler order}
(or, in short, \emph{Rudin-Keisler compatible})
if there is an increasing $h\in\NN$ such that $\cF_1/h\cup \cF_2/h$
satisfies the finite intersection property (that is, it is a filter base).
If $\cF_1,\cF_2$ are Rudin-Keisler compatible ultrafilters, then there
is an increasing $h\in\NN$ such that $\cF_1/h=\cF_2/h$.

\begin{defn}
\emph{NCF (near coherence of filters)} is the assertion that every two nonprincipal
ultrafilters on $\N$ are Rudin-Keisler compatible.
\end{defn}
NCF is independent of ZFC \cite{ShBl, NCFIII}, and has many equivalent forms and
implications (e.g., \cite{NCFI, NCFII}).

In the sequel, we often use the following convenient notation
for $f,g\in\NN$:
$$[f\le g] = \{n : f(n)\le g(n)\}.$$

\begin{thm}[Bartoszy\'nski-Shelah-Tsaban \cite{huremen2}]\label{dfin-finite}
NCF holds if, and only if, $\Dfin$ is additive.
\end{thm}
\begin{proof}
($\Impl$)
Assume that $Y_1,Y_2\in\Dfin$.
We may assume that all elements of $Y_1$ and $Y_2$ are strictly increasing
and that $Y_1$ and $Y_2$ are closed under finite maxima.
Thus, it suffices to show that
$$\{\max\{f_1,f_2\} : f_1\in Y_1,\ f_2\in Y_2\}$$
is not dominating.
For each $i=1,2$, do the following: Choose an increasing $g_i\in\NN$ witnessing that $Y_i$ is not
dominating. The set $\{[f\le g] : f\in Y_i\}$ has the finite intersection property.
Extend it to a nonprincipal ultrafilter $\cF_i$.

Fix an increasing $h\in\NN$ such that $\cF_1/h\cup \cF_2/h$
has the finite intersection property. Define $g\in\NN$ by
$g(n)=\max\{g_1(h(n+1)),g_2(h(n+1))\}$ for each $n$.
Given $f_1\in Y_1, f_2\in Y_2$, let $a$ be the infinite set
$[f_1\le g_1]/h\cap [f_2\le g_2]/h$. For each $n\in a$
and each $i=1,2$, there is $k\in\intvl{h}{n}$ such that
$f_i(k)\le g_i(k)$. Thus,
$$f_i(n)\le f_i(h(n))\le f_i(k)\le g_i(k)\le g_i(h(n+1))\le g(n),$$
thus $\max\{f_1(n),f_2(n)\}\le g(n)$ for all $n\in a$.

\medskip

($\Leftarrow$) We will use the following.

\begin{lem}[\cite{huremen2}]\label{distantdisjoint}
If NCF fails, then there exist ultrafilters $\cF_1$ and $\cF_2$
such that for each increasing $h\in\NN$ there exist
$a_1\in\cF_1/h$ and $a_2\in\cF_2/h$ such that for all
$n\in a_1$ and $m\in a_2$, $|n-m|>1$.
\end{lem}
\begin{proof}
Assume that $\cF_1$ and $\cF_2$ are Rudin-Keisler \emph{in}compatible nonprincipal ultrafilters
and let $h$ be an increasing element of $\NN$.
Define increasing $f_0,f_1\in\NN$ by
\begin{eqnarray*}
f_0(n) & = & h(2n)\\
f_1(n) & = & h(2n+1)
\end{eqnarray*}
Then there exist
\begin{eqnarray*}
X_1\in \cF_1/f_0 &~& X_2\in\cF_2/f_0\\
Y_1\in \cF_1/f_1 &~& Y_2\in\cF_2/f_1
\end{eqnarray*}
such that the sets $X_1\cap X_2=Y_1\cap Y_2=\emptyset$.\footnote{Since
nonprincipal filters are closed under finite modifications,
we can shrink the elements to turn the finite intersection into an empty
intersection.
}
For $i=1,2$ let
\begin{eqnarray*}
\tilde X_i & = & 2\cdot X_i\cup (2\cdot X_i+1)\\
\tilde Y_i & = & (2\cdot Y_i+1)\cup (2\cdot Y_i+2)
\end{eqnarray*}
Observe that $\tilde X_1\cap \tilde X_2=\tilde Y_1\cap \tilde Y_2=\emptyset$ either.
Now,
\begin{eqnarray*}
\Union_{n\in X_i}\intvl{f_0}{n} & = & \Union_{n\in \tilde X_i}\intvl{h}{n}\\
\Union_{n\in Y_i}\intvl{f_1}{n} & = & \Union_{n\in \tilde Y_i}\intvl{h}{n}
\end{eqnarray*}
therefore $\tilde X_i,\tilde Y_i\in\cF_i/h$, thus
$a_i = \tilde X_i\cap \tilde Y_i\in\cF_i/h$.
If $n\in a_1$ is even, then $n,n+1\in\tilde X_1$, and $n-1,n\in\tilde Y_1$.
Thus, if $n$ is large enough, then  $n,n+1\nin\tilde X_2$, and $n-1,n\nin\tilde Y_2$,
therefore $n-1,n,n+1\nin a_2$.
The case that $n\in a_1$ is odd is similar.
\end{proof}

For a filter $\cF$ and an increasing $g\in\NN$, define
$$Y_{\cF,g}=\{f\in\NN : [f\le g]\in\cF\}.$$
Then $Y_{\cF,g}\in\Dfin$.
It therefore suffices to prove the following.

\begin{lem}\label{addDfinD2}
If $\cF_1$ and $\cF_2$ are as in Lemma \ref{distantdisjoint},
and $g(n)\ge 2n$ for each $n$,
then $Y_{\cF_1,g}\cup Y_{\cF_2,g}$ is $2$-dominating.
\end{lem}
\begin{proof}
Let $f\in\NN$ be any increasing function.
Define by induction
\begin{eqnarray*}
h(0) & = & 0\\
h(n+1) & = & f(h(n))+1
\end{eqnarray*}
By the assumption,
there exist $a_1\in\cF_1/h$ and $a_2\in\cF_2/h$ such that
for each $n\in a_1$ and $m\in a_2$, $|n-m|>1$.

Fix $i<2$.
For each $n$, define
$$f_i(n) =
\begin{cases}
f(h(k-1))+n-h(k-1) & n\in \intvl{h}{k} \mbox{ for }k\in a_i\\
f(h(k))+n-h(k)     & \!\!\begin{array}{l}n\in \intvl{h}{k}\\ \mbox{ where }k\nin a_i, k+1\in a_i\end{array}\\
f(n) & \mbox{otherwise}
\end{cases}
$$
It is not difficult to verify that $f_i$ is increasing.

For each $k\in a_i$ and $n\in \intvl{h}{k}$,
\begin{eqnarray*}
f_i(n) & = & f(h(k-1))+n-h(k-1) \le\\
& \le & h(k)+n-h(k-1)\le h(k)+n\le 2n\le g(n).
\end{eqnarray*}
Therefore $f_i\in Y_{\cF_i,g}$.

For each $n$ let $k$ be such that $n\in\intvl{h}{k}$.
If $n$ is large enough, then either $k,k+1\nin a_1$,
and therefore $f_1(n) = f(n)$, or else
$k,k+1\nin a_2$, and therefore $f_2(n) = f(n)$, that is,
$f(n)\le\max\{f_1(n),f_2(n)\}$.\footnote{In fact we get equality here.}
\end{proof}
This completes the proof of Theorem \ref{dfin-finite}.
\end{proof}

Let $\add(\Dfin,\fD)$ denote the minimal cardinality of
a collection of members of $\Dfin$ whose union is dominating.
It is immediate that $\fb\le\add(\Dfin,\fD)$.

\begin{lem}[Blass \cite{BlassNew}]\label{Blassbg}
$\max\{\fb,\fg\}\le\add(\Dfin,\fD)$.
\end{lem}
\begin{proof}
We need only prove that $\fg\le\add(\Dfin,\fD)$.
Assume that $\kappa<\fg$ and $Y_\alpha\in\Dfin$, $\alpha<\kappa$.
We may assume each $Y_\alpha$ is closed under
pointwise maxima of its finite subsets.
For each $\alpha$, let $g_\alpha$ be a witness for $Y_\alpha$ not being
dominating, and extend $\{[f\le g_\alpha] : f\in Y_\alpha\}$ to a nonprincipal
ultrafilter $\cF_\alpha$ on $\N$.

We will use the following ``morphism''.

\begin{lem}[Mildenberger \cite{Mild01, MShT:847}]\label{morph}
For each $f\in\NN$ and each ultrafilter $\cU$,
$$\cG_{\cU,f}=\{a\in\roth : f\le_\cU\Next{a}\}$$
is groupwise dense.
\end{lem}
\begin{proof}
Clearly, $\cG_{\cU,f}$ is closed under taking almost subsets.
Assume that $\{\intvl{h}{n}: n\in\w\}$ is an interval partition of $\w$.
By merging consecutive intervals we may assume that for each $n$,
and each $k\in\intvl{h}{n}$, $f(k)\le h(n+2)$.

Since $\cU$ is an ultrafilter, there exists $\ell\in\{0,1,2\}$ such that
$$a_\ell=\Union_n \intvl{h}{3n+\ell}\in\cU$$
Take $a=a_{\ell+2\bmod 3}$.
For each $k\in a_\ell$, let $n$ be such that $k\in [h(3n+\ell),\allowbreak h(3n+\ell+1))$.
Then $f(k)\le h(3n+\ell+2)=\Next{a}(k)$. Thus $a\in\cG_{\cU,f}$.
\end{proof}
Thus, we can take $a\in\bigcap_{\alpha<\kappa}\cG_{\cU_\alpha,g_\alpha}$,
and $g=\Next{a}$ will witness that $\Union_{\alpha<\kappa}Y_\alpha$ is not
dominating.
\end{proof}

\begin{thm}\label{dfin-many}
If $\Dfin$ is additive (equivalently, NCF holds),
then it is $\add(\Dfin,\fD)$-additive and therefore $\max\{\fb,\fg\}$-additive.
In particular, in this case it is $\sigma$-additive.
\end{thm}
\begin{proof}
Assume that $\kappa<\add(\Dfin,\fD)$ and $Y_\alpha\in\Dfin$, $\alpha<\kappa$.
We may assume that each $Y_\alpha$ is closed under pointwise maxima of finite subsets,
and that the family $\{Y_\alpha : \alpha<\kappa\}$ is additive. It follows that
$$\maxfin\left(\Union_{\alpha<\kappa} Y_\alpha\right ) = \Union_{\alpha<\kappa} Y_\alpha$$
and is therefore not dominating. Thus, $\Union_{\alpha<\kappa} Y_\alpha\in\Dfin$.

The second assertion follows from Lemma \ref{Blassbg}.
\end{proof}

\begin{cor}\label{estMen}
If NCF holds, then
$$\max\{\fb,\fg\}\le\add(\ufin(\cO,\Omega))
\le\add(\sone(\BG,\BO))\le\cf(\fd)=\fd.\qed$$
\end{cor}

Recently, Banakh and Zdomskyy improved Theorem \ref{dfin-many} and Corollary \ref{estMen},
by showing that NCF implies that $\add(\Dfin)=\add(\ufin(\cO,\Omega))=\fd$.

\begin{prob}
Is any of the classes
$\sfin(\Omega,\Omega)$, $\sone(\Gamma,\Omega)$, and $\sfin(\Gamma,\allowbreak\Omega)$
consistently additive?
\end{prob}

For the Borel case there remains only one unsolved class.
\begin{prob}
Is $\sfin(\BO,\BO)$ consistently additive?
\end{prob}

\subsection{On splitting properties}

\begin{thm}[Zsomskyy \cite{SF1, combimgs}]
It is consistent that $\add(\Split(\Lambda,\allowbreak\Lambda))=\add(\Split(\BL,\BL))=\fb=\fu$.
\end{thm}
\begin{proof}
In \cite{SF1, combimgs} it is proved that
$\fu<\fg$ implies that $\Split(\Lambda,\Lambda)=\ufin(\cO,\Gamma)$, and
the same assertion holds in the Borel case.
The theorem follows from Corollary \ref{adds}, together with the
fact that $\fu<\fg$ implies that $\fb=\fu$ \cite{BlassHBK}.
\end{proof}

The last theorem implies that one cannot obtain a negative solution to
Problem \ref{SpLamLamAdd} in ZFC.

\section{$\tau$-covers}
$\cU$ is a \emph{$\tau$-cover} of $X$ if it is a large cover of $X$
(that is, each member of $X$ is contained in infinitely many members of the cover),
and for all $x,y\in X$, (at least)
one of the sets $\{U\in\cU : x\in U, y\nin U\}$ and
$\{U\in\cU : y\in U, x\nin U\}$ is finite.
$\tau$-covers are motivated by the \emph{tower number} $\ft$
\cite{tau} and were incorporated into the framework of selection principles
in \cite{tautau}.
Every open $\tau$-cover of a set of reals contains a countable $\tau$-cover of
that set \cite{split}.
Let $\Tau$ and $\BT$
denote the collections of countable open and Borel $\tau$-covers of $X$, respectively.

\subsection{On the Scheepers diagram}

\begin{figure}[!ht]
{\tiny
\begin{changemargin}{-3cm}{-3cm}
\begin{center}
$\xymatrix@C=7pt@R=6pt{
&
&
& \sr{\ufin(\cO,\Gamma)}{\fb}\ar[r]
& \sr{\ufin(\cO,\Tau)}{{\max\{\fb,\fs\}}}\ar[rr]
&
& \sr{\ufin(\cO,\Omega)}{\fd}\ar[rrrr]
&
&
&
& \sr{\sfin(\cO,\cO)}{\fd}
\\
&
&
& \sr{\sfin(\Gamma,\Tau)}{\fb}\ar[rr]\ar[ur]
&
& \sr{\sfin(\Gamma,\Omega)}{\fd}\ar[ur]
\\
\sr{\sone(\Gamma,\Gamma)}{\fb}\ar[uurrr]\ar[rr]
&
& \sr{\sone(\Gamma,\Tau)}{\fb}\ar[ur]\ar[rr]
&
& \sr{\sone(\Gamma,\Omega)}{\fd}\ar[ur]\ar[rr]
&
& \sr{\sone(\Gamma,\cO)}{\fd}\ar[uurrrr]
\\
&
&
& \sr{\sfin(\Tau,\Tau)}{{\min\{\fb,\fs\}}}\ar'[r][rr]\ar'[u][uu]
&
& \sr{\sfin(\Tau,\Omega)}{\fd}\ar'[u][uu]
\\
\sr{\sone(\Tau,\Gamma)}{\ft}\ar[rr]\ar[uu]
&
& \sr{\sone(\Tau,\Tau)}{\ft}\ar[uu]\ar[ur]\ar[rr]
&
& \sr{\sone(\Tau,\Omega)}{\fod}\ar[uu]\ar[ur]\ar[rr]
&
& \sr{\sone(\Tau,\cO)}{\fod}\ar[uu]
\\
&
&
& \sr{\sfin(\Omega,\Tau)}{\fp}\ar'[u][uu]\ar'[r][rr]
&
& \sr{\sfin(\Omega,\Omega)}{\fd}\ar'[u][uu]
\\
\sr{\sone(\Omega,\Gamma)}{\fp}\ar[uu]\ar[rr]
&
& \sr{\sone(\Omega,\Tau)}{\fp}\ar[uu]\ar[ur]\ar[rr]
&
& \sr{\sone(\Omega,\Omega)}{\cov(\cM)}\ar[uu]\ar[ur]\ar[rr]
&
& \sr{\sone(\cO,\cO)}{\cov(\cM)}\ar[uu]
}$
\end{center}
\end{changemargin}
}

\caption{The Scheepers diagram, enhanced with $\tau$-covers}\label{tauSch}
\end{figure}
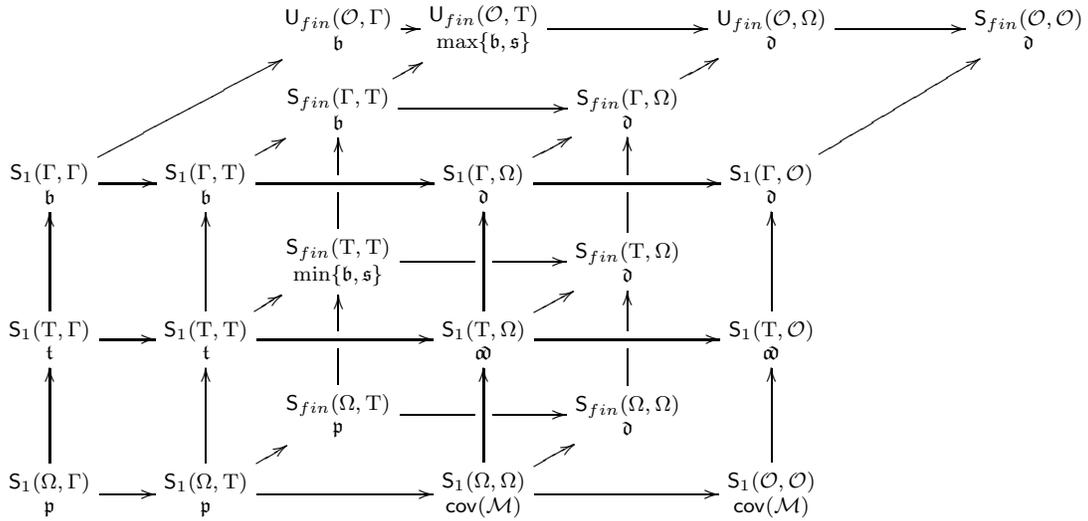

Taking $\Tau$ into account and removing trivial properties and known equivalences,
we obtain the diagram in Figure \ref{tauSch} \cite{tautau, MShT:858}.
In this diagram too, the critical cardinality of each property appears below it.
A similar diagram, with several additional equivalences, is available in the
Borel case \cite{tautau}.

\begin{prop}
$\sone(\Tau,\cO)$ and $\sone(\BT,\B)$ are $\sigma$-additive.
\end{prop}
\begin{proof}
As in Proposition \ref{easyadd}.
\end{proof}

\begin{defn}
For each countable cover of $X$ enumerated bijectively as
$\cU=\seq{U_n}$ we associate the \emph{Marczewski function}
$h_\cU: X\to P(\N)$, defined by $h_\cU(x) = \{ n : x\in U_n\}$
for each $x\in X$.
\end{defn}

$\cU$ is a large cover of $X$ if, and only if, $h_\cU[X]\sbst\roth$.
Let $Y\sbst\roth$.
$Y$ is \emph{centered} if for each finite $F\sbst Y$,
$\cap F$ is infinite. A set $a\in\roth$ is a \emph{pseudo-intersection}
of $Y$ if $a\as y$ for all $y\in Y$.
$Y$ is \emph{linearly quasiordered by $\as$} if for all $y,z\in Y$,
$y\as z$ or $z\as y$. Note that if $Y$ has a pseudo-intersection or
is linearly quasiordered by $\as$, then $Y$ is centered.

\begin{lem}[Tsaban \cite{tau}]\label{notions}
Assume that $\cU$ is a countable large cover of $X$.
\be
\itm $\cU$ is an $\w$-cover of $X$ if, and only if, $h_\cU[X]$ is centered.
\itm $\cU$ contains a $\gamma$-cover of $X$ if, and only if,
$h_\cU[X]$ has a pseudo-intersection.
\itm $\cU$ is a $\tau$-cover of $X$ if, and only if, $h_\cU[X]$ is linearly quasiordered by $\as$.\hfill\qed
\ee
\end{lem}

For families $\scrB\sbst\scrA$ of covers of a space $X$, define the property
\emph{$\scrA$ choose $\scrB$} as follows.
\bi
\item[$\binom{\scrA}{\scrB}$:]
For each $\cU\in\scrA$, there is $\cV\sbst\cU$ such that $\cV\in\scrB$.
\ei
This is a prototype for many classical topological notions, most notably compactness and being Lindel\"of.

\begin{thm}[Tsaban \cite{tautau}]\label{addT}
$\add(\binom{\Tau}{\Gamma})=\add(\binom{\BT}{\BG})=\ft$.
\end{thm}
\begin{proof}
We prove the open case.
Assume that $\kappa<\ft$,
and let $X_\alpha$, $\alpha<\kappa$, be sets satisfying
$\binom{\Tau}{\Gamma}$. Let $\cU$ be a countable open $\tau$-cover
of $X=\Union_{\alpha<\kappa}X_\alpha$.
By Lemma \ref{notions},
$h_\cU[X]=\Union_{\alpha<\kappa}h_\cU[X_\alpha]$ is linearly
quasiordered by $\as$. Since each $X_\alpha$ satisfies
$\binom{\Tau}{\Gamma}$, for each $\alpha$ $\cU$ contains a
$\gamma$-cover of $X_\alpha$, that is, $h_\cU[X_\alpha]$ has a
pseudo-intersection.

\begin{lem}[Tsaban \cite{tau}]\label{p.i.pieces}
Assume that $Y\sbst\roth$ is linearly quasiordered by $\as$, and
for some $\kappa<\ft$, $Y=\Union_{\alpha<\kappa}Y_\alpha$ where each
$Y_\alpha$ has a pseudo-intersection.
Then $Y$ has a pseudo-intersection.
\end{lem}
\begin{proof}
If for each $\alpha<\kappa$ $Y_\alpha$ has a pseudo-intersection $y_\alpha\in Y$,
then a pseudo-intersection of $\{y_\alpha : \alpha<\kappa\}$
is also a pseudo-intersection of $Y$.
Otherwise, there exists $\alpha<\kappa$ such that $Y_\alpha$ has no pseudo-intersection
$y\in Y$.
That is, for all $y\in Y$ there exists a $z\in Y_\alpha$ such that
$y\not\as z$; thus
$z\as y$.
Therefore, a pseudo-intersection of $Y_\alpha$ is also a pseudo-intersection of $Y$.
\end{proof}
By Lemma \ref{p.i.pieces}, $h_\cU[X]$ has a
pseudo-intersection, that is, $\cU$ contains a $\gamma$-cover of
$X$.
\end{proof}

\forget
The following \emph{cancellation laws} are useful.
\begin{lem}[Tsaban \cite{tautau}]\label{cancellation}
For collections of covers $\scrA\spst\scrB\spst\scrC$:
\be
\itm $\binom{\scrA}{\scrB}\cap\binom{\scrB}{\scrC}=\binom{\scrA}{\scrC}$.
\itm $\binom{\scrA}{\scrB}\cap\sfin(\scrB,\scrC)=\sfin(\scrA,\scrC)$.
\itm $\sfin(\scrA,\scrB)\cap\binom{\scrB}{\scrC}=\sfin(\scrA,\scrC)$.
\itm $\binom{\scrA}{\scrB}\cap\sone(\scrB,\scrC)=\sone(\scrA,\scrC)$.
\itm If $\scrC$ is closed under taking supersets, then
$\sone(\scrA,\scrB)\cap\binom{\scrB}{\scrC}=\sone(\scrA,\scrC)$.\hfill\qed
\ee
\end{lem}
\forgotten

\begin{cor}\label{adds1tg}
$\add(\sone(\Tau,\Gamma))=\add(\sone(\BT,\BG))=\ft$.
\end{cor}
\begin{proof}
Note that
$$\sone(\Tau,\Gamma) = \binom{\Tau}{\Gamma}\cap\sone(\Gamma,\Gamma).$$
It follows that $\add(\sone(\Tau,\Gamma))$ is at least the minimum of
the additivity numbers of $\binom{\Tau}{\Gamma}$ and $\sone(\Gamma,\Gamma)$,
which are $\ft$ (Theorem \ref{addT}) and $\fh$ (Theorem \ref{s1gg}), respectively.
As $\ft\le\fh$ \cite{BlassHBK}, $\add(\sone(\Tau,\Gamma))\ge\ft$.
On the other hand, $\add(\sone(\Tau,\Gamma))\le\non(\sone(\Tau,\Gamma))=\ft$
(Figure \ref{tauSch}).

In the Borel case use $\add(\sone(\BG,\BG))=\fb\ge\ft$ (Theorem \ref{adds}).
\end{proof}

Note that $\sfin(\Omega,\Tau)$ implies $\binom{\Omega}{\Tau}$.

\begin{cor}[Tsaban \cite{tautau}]\label{notadd2}
Assume \CH{}. Th\-en no class between $\sone(\BO,\BG)$ and $\binom{\Omega}{\Tau}$ (inclusive)
is additive.
\end{cor}
\begin{proof}
By Theorem \ref{gammanotadd}, there are sets $A$ and $B$ satisfying $\sone(\BO,\BG)$, such that
$A\cup B$ does not satisfy $\sone(\Omega,\Gamma)$.
Now,
$$\sone(\Omega,\Gamma)=\binom{\Omega}{\Tau}\cap\sone(\Tau,\Gamma),$$
and by Corollary \ref{adds1tg}, $A\cup B$ satisfies $\sone(\Tau,\Gamma)$.
Thus, $A\cup B$ does not satisfy $\binom{\Omega}{\Tau}$.
\end{proof}

\begin{prob}
Is any of the properties
$\sone(\Tau,\Tau)$, $\sfin(\Tau,\Tau)$,
$\sone(\Gamma,\Tau)$, $\sfin(\Gamma,\Tau)$,
and $\ufin(\cO,\Tau)$ (or any of their Borel
versions) provably (or at least consistently)
additive?
\end{prob}

Zdomskyy \cite{SF2} proved that consistently, $\ufin(\cO,\Tau)=\ufin(\cO,\Gamma)$, and in particular,
$\ufin(\cO,\Tau)$ is consistently additive. Mildenberger, Shelah, and Tsaban \cite{MShT:858}
proved that $\sone(\Tau,\Tau)$ is additive if, and only if, $\sone(\Tau,\Tau)=\sone(\Tau,\Gamma)$.
We do not know whether the latter assertion is consistent.

\begin{prob}
Is any of the classes
$\sfin(\Omega,\Tau)$,
$\sone(\Tau,\Omega)$, and
$\sfin(\Tau,\Omega)$
consistently additive?
\end{prob}

\subsection{On splitting properties}

Here, taking $\Tau$ into account and removing trivialities and
equivalences, we obtain the following diagram (in the open case,
and a similar one in the Borel case) \cite{split}:
$$\xymatrix{
\sr{\Split(\Lambda, \Lambda)}{\fr} \ar[r] & \sr{\Split(\Omega, \Lambda)}{\fu} \ar[r] & \sr{\Split(\Tau, \Tau)}{\mbox{undefined}}\\
                           & \sr{\Split(\Omega, \Omega)}{\fu}\ar[u]\\
                          & \sr{\Split(\Omega, \Tau)}{\fp}\ar[u]\\
\sr{\Split(\Omega, \Gamma)}{\fp} \ar[uuu]\ar[ur]\ar[rr]     & & \sr{\Split(\Tau,\Gamma)}{\ft}\ar[uuu]\\
}$$
We also have that $\Split(\Tau,\Gamma)=\binom{\Tau}{\Gamma}$ \cite{split}.
By Theorem \ref{addT}, $\add(\Split(\Tau,\allowbreak\Gamma))=\ft$.

\begin{thm}[Tsaban \cite{split}]\label{utt}
$\fu\le\add(\Split(\Tau,\Tau))$. 
\end{thm}
\begin{proof}
A nonprincipal ultrafilter $U$ on $\N$ is called a \emph{simple $P$-point}
if there exists a base $B$ for $U$ such that $B$ is linearly quasiordered
by $\as$. We call such a base a \emph{simple $P$-point base}.
\begin{lem}[\cite{split}]\label{Ppoint}
$X$ satisfies $\Split(\Tau,\Tau)$ if, and only if,
for each countable open $\tau$-cover $\cU$ of $X$,
$h_\cU[X]$ is not a simple $P$-point base.\hfill\qed
\end{lem}
Thus, our theorem follows from the following Ramseyan property.
\begin{lem}[\cite{split}]
Assume that $\lambda<\fu$ and $B = \Union_{\alpha<\lambda}B_\alpha$
is a simple $P$-point base. Then there exists $\alpha<\lambda$
such that $B_\alpha$ is a simple $P$-point base.
\end{lem}
\begin{proof}
Assume that $B$ is a simple $P$-point base
and $U$ is the simple $P$-point it generates.
In particular, $B$ is linearly ordered by $\as$.
We will show that some $B_\alpha$ is a base for $U$.
Assume otherwise.
For each $\alpha<\lambda$ choose $a_\alpha\in U$ that
witnesses that $B_\alpha$ is not a base for $U$,
and $\tilde a_\alpha\in B$ such that $\tilde a_\alpha\as a_\alpha$.
As $B$ is linearly ordered by $\as$,
$\tilde a_\alpha$ is a pseudo-intersection of $B_\alpha$.

The cardinality of the linearly ordered set $Y = \{\tilde a_\alpha : \alpha<\lambda\}$
is smaller than $\fu$. Thus it is not a base for $U$ and we can find again
an element $a\in\cF$ which is a pseudo-intersection of $Y$, and therefore of $B$;
a contradiction.
\end{proof}
This completes the proof of Theorem \ref{utt}.
\end{proof}

Consistently, there are no $P$-points \cite{barju}.
By Lemma \ref{Ppoint}, in such a model $\Split(\Tau,\Tau)=P(\R)$ and therefore
$\add(\Split(\Tau,\Tau))$ is undefined.

Note that $\Split(\Omega, \Tau)$ implies $\binom{\Omega}{\Tau}$,
and since $\Split(\BO, \BG)=\binom{\BO}{\BG}=\sone(\BO,\BG)$,
we have by Corollary \ref{notadd2} that no class between $\Split(\BO, \BG)$ and
$\Split(\Omega, \Tau)$ (inclusive) is provably additive.

Thus, $\Split(\Omega, \Lambda)$, $\Split(\Tau, \Tau)$, and $\Split(\Tau,\Gamma)$
are (provably) $\sigma$-addi\-tive, whereas
$\Split(\Omega, \Omega)$, $\Split(\Omega, \Tau)$, and $\Split(\Omega, \Gamma)$ are not provably additive.
The situation for $\Split(\Lambda, \Lambda)$ is Problem \ref{SpLamLamAdd}.

\begin{prob}
Improve the lower bound or the upper bound in the inequality
$\aleph_1\le\add(\Split(\Omega, \Lambda))\le\fc$.
\end{prob}

\begin{prob}
Can the lower bound $\fu$ on $\add(\Split(\Tau, \Tau))$ be improved?
\end{prob}

\subsection*{Acknowledgements}
We thank Assaf Rinot and Lyubomyr Zdomskyy for reading the paper and making useful comments.

\end{document}